\documentclass[12pt]{article}
\usepackage{graphicx}
\usepackage{amsfonts}
\usepackage{amssymb}

\input xy
\xyoption{all}

\newtheorem{theorem}{Theorem}[section]
\newtheorem{lemma}[theorem]{Lemma}
\newtheorem{example}[theorem]{Example}
\newtheorem{proposition}[theorem]{Proposition}

\newcommand{\Ker}{\mbox{\rm Ker\,}}
\newcommand{\Coker}{\mbox{\rm Coker\,}}
\newcommand{\Ima}{\mbox{\rm Im\,}}
\newcommand{\End}{\mbox{\rm End}}
\newcommand{\Ann}{\mbox{\rm Ann\,}}
\newcommand{\mod}{\mbox{\rm mod\,}}

\newcommand{\mA}{\mbox{$\mathbb A$}}
\newcommand{\mD}{\mbox{$\mathbb D$}}
\newcommand{\mN}{\mbox{$\mathbb N$}}

\newcommand{\udim}{\underline{\dim}\,}
\newcommand{\ul}[1]{\underline{#1}}

\newcommand{\vertex}{\circle*{3}}
\newcommand{\bvertex}{\circle*{5}}

\newcommand{\bwertex}{\circle{5}}
\newcommand{\svertex}{\circle*{1}}

\newcommand{\Hom}[1]{\textrm{Hom}_{#1}}
\newcommand{\Ext}[2]{{\rm Ext}^{#1}_{#2}}

\begin{document}
\title{Gentle algebras arising from surface triangulations}
\author{Ibrahim Assem \and Thomas Br\"ustle \and Gabrielle
Charbonneau-Jodoin \and Pierre-Guy Plamondon\thanks{Ibrahim Assem is
partially supported by NSERC of Canada and the Universit\'e de
Sherbrooke.  Gabrielle Charbonneau  was working under a summer research
fellowship of NSERC. Thomas Br\"ustle is partially supported by NSERC,
by Bishop's University and the Universit\'e de Sherbrooke.  Pierre-Guy
Plamondon was supported by an NSERC graduate fellowship}}
  \date{}
\maketitle

\begin{abstract}
In this paper, we associate an algebra $A(\Gamma)$ to a triangulation $\Gamma$ of a surface $S$ with a set of boundary marking points. This algebra $A(\Gamma)$ is gentle and Gorenstein of dimension one. We also prove that $A(\Gamma)$ is cluster-tilted if and only if it is cluster-tilted of type $\mA$ or $\widetilde{\mA}$, or if and only if the surface $S$ is a disc or an annulus. Moreover all cluster-tilted algebras of type $\mA$ or $\widetilde{\mA}$ are obtained in this way.

\end{abstract}

\section{Introduction}
Among the main recent results in the fast-growing theory of cluster algebras is the paper of Fomin, Shapiro and Thurston \cite{FST}, relating triangulations of oriented surfaces to cluster algebras. 
This approach, which existed since the beginning of the theory \cite{CCS}, was followed, among others, in \cite{La,Sch}.
In the same spirit, we consider in the present paper an unpunctured oriented surface $S$, and a finite set of points $M$ lying on the boundary of $S$, and intersecting every boundary component of $S$.
We then associate to a triangulation $\Gamma$ of the marked surface $(S,M)$ a quiver $Q(\Gamma)$, and a potential on $Q(\Gamma)$ (in the sense of \cite{DWZ}), thus defining an algebra $A(\Gamma)$, namely the (non-completed) Jacobian algebra defined by $Q(\Gamma)$ and the associated potential.
\bigskip

Such an algebra $A(\Gamma)$ has some very nice properties: it is always Gorenstein of dimension one, and also it is a gentle algebra in the sense of \cite{AS}. 
In the unpunctured case studied here, our definition coincides with Labardini's definition of a quiver with potential associated to a (possibly punctured) surface  \cite{La}.
But in the punctured case, one does not get gentle algebras, or even string algebras.
For instance, a once-punctured disc gives rise to cluster-tilted algebras of type $\mD$, see \cite{Sch}.
\medskip

Gentle algebras form a particularly nice subclass of the class of string algebras of \cite{BR} and are much investigated in the representation theory of algebras.
For instance, this subclass contains the tilted algebras of type $\mA$ and $\widetilde{\mA}$ (see \cite{A} and \cite{AS}, respectively) and it is closed under tilting and even under derived equivalence (see \cite{S} and \cite{SZ}, respectively).
\smallskip

Our objective in this paper is twofold. Firstly, we ask which gentle algebras arise in this way, that is, are induced from triangulations of an unpunctured surface with boundary marked points. 
We show in (2.8) below that this is the case for every gentle algebra such that every relation lies on what we call a 3-cycle with radical-square zero (see (2.2) or \cite{BV} for the definition).
Secondly, we ask which gentle algebras are cluster-tilted. 
The class of cluster-tilted algebras, introduced in \cite{BMR}, has been much investigated and is by now well-understood (see, for instance, \cite{ABS,BFPPT,BMRRT,BV,CCS,Ke2,KR,Sch}).
In particular, it was shown in \cite{ABS} that every cluster-tilted algebra is the relation-extension of a tilted algebra, that is, it is the trivial extension of a tilted algebra $C$ by the $C-C-$bimodule $\Ext{2}{C}(DC,C)$. We may now state the main result of this paper.

\begin{theorem}\label{main}
Let $A(\Gamma)$ be the algebra associated to the triangulation $\Gamma$ of an unpunctured marked  surface
$(S,M)$.  Then the following statements are equivalent:

\begin{itemize}
\item[(1)] $A(\Gamma)$ is cluster-tilted
\item[(2)] $A(\Gamma)$ is cluster-tilted of type $\mA$ and $\widetilde{\mA}$
\item[(3)]  $A(\Gamma)$ is the relation-extension of a tilted algebra of type $\mA$ and $\widetilde{\mA}$
\item[(4)] the surface $S$ is a disc or an annulus.
\end{itemize}
Moreover, all cluster-tilted algebras of type  $\mA$ (or $\widetilde{\mA}$) are of the form $A(\Gamma)$ for some triangulation of a disc $S$ (or an annulus $S$, respectively).
\end{theorem}

Actually, we prove in (3.3) that a cluster-tilted algebra is gentle if and only if it is of type $\mA$ and $\widetilde{\mA}$, or if and only if it is the relation-extension of a gentle tilted algebra, and the latter coincide with the tilted algebras of type $\mA$ or $\widetilde{\mA}$, respectively. 
\medskip

The case where $S$ is a disc has already been studied in \cite{CCS}, and it is known that the bound quivers of all cluster-tilted algebras of type $\mA$ arise from triangulations of the (unpunctured) disc. 
These algebras have also been described explicitly in \cite{BV}. Also, the potential we use for defining the cluster-tilted algebras of type $\widetilde{\mA}$ is a particular case of the potential recently defined by B. Keller \cite{Ke2}. However, we do not use this fact, but rather present another proof (anterior to Keller's result), which uses \cite{ABS} and properties of the second extension group.
\bigskip

The paper is organised as follows: in section 2, we define our algebras $A(\Gamma)$ and prove their main properties in (2.7). Section 3 is devoted to the classification of the gentle cluster-titled algebras and section 4 to the proof of our main theorem and some of its consequences. We end the paper with an example of an algebra $A(\Gamma)$ which is not of polynomial growth (in the sense of \cite{Sk}).


\section{Algebras arising from surface triangulations}

Throughout this paper, the algebras we consider are basic  connected
algebras over a fixed algebraically closed field $k$. Unless otherwise
stated, all algebras are finite-dimensional. Consequently, they are
given in the form  $A = kQ/I$ where $Q$ is a quiver and $I$ is an
admissible ideal of the path algebra $kQ$, see \cite{ASS}.
The pair $(Q,I)$  is called a {\em bound quiver}, and the algebra $A = kQ/I$
is referred to as a {\em bound quiver algebra}.

Given  a bound quiver algebra $A=kQ/I$, for every vertex $x$ of $Q$ we
denote by $e_x$ the idempotent of $A$ associated to $x$.  Also, $P_x$,
$I_x$ and $S_x$ will be the corresponding indecomposable projective
module, indecomposable injective module and simple module, respectively.
\bigskip

We study in this section the  algebra  associated with a surface
triangulation. For background material on oriented surfaces we refer to
\cite{Ma}.

\subsection{The medial quiver $Q(\Gamma)$}
We first recall from \cite{FST} the construction of a quiver  for every
triangulation of a marked surface: Let $S$ be an oriented surface with
boundary $\partial S$, and let $M$ be a non-empty finite set of points on
$\partial S$ intersecting each connected component of  the boundary
$\partial S$. In this paper, we only consider the case where there are no punctures, that
is we request that the set of marked points $M$ be contained in the
boundary $\partial S$.
The pair $(S,M)$ is referred to as an {\em unpunctured  bordered surface with
marked points.}

An {\em arc}  in $(S, M)$ is a curve $\gamma$ in $S$ such that
\begin{itemize}
\item the endpoints of $\gamma$ are marked points in $M$
\item $\gamma$ does not intersect itself, except that its endpoints may
coincide
\item $\gamma$ intersects  the boundary of $\partial S$ only in its
endpoints
\item $\gamma$ does not cut out a monogon  (that is,  $\gamma$ is not
contractible into a point of $M$).
\end{itemize}

We call an arc $\gamma$ a {\em boundary arc} if it cuts out a  digon (that is,
 $\gamma$ is homotopic to a curve $\delta$ on the boundary $\partial S$
 that intersects $M$ only in its endpoints).
Otherwise, $\gamma$ is said to be an {\em internal arc.}
Each arc $\gamma$ is considered up to homotopy in the class of such curves.
A {\em triangulation} of $(S,M)$ is a maximal collection $\Gamma$ of arcs that
do not intersect in the interior of $S$ (more precisely, there are
curves in their respective homotopy classes that  do not intersect in the
interior of $S$).

\begin{proposition}[{\cite[(2.10)]{FST}}]  In each  triangulation of
$(S,M)$, the number of internal arcs is
$$n = 6g + 3b  + c - 6 $$
where $g$ is the genus of $S$, $b$ is the number of boundary components,
and $c=|M|$ is the number of marked points.
\end{proposition}

This Proposition also indicates that in  some cases  a triangulation does
not exist (for instance a disc with one marked point would give $n=
-2$). We consider from now on only marked surfaces $(S,M)$ that admit  a
triangulation.
Given a triangulation $\Gamma$, we also refer to $M$ as the set of
vertices of  $\Gamma$. The triangles are the components of $S \backslash
\Gamma$ with the arcs of $\Gamma$ as edges.

We denote by $Q(\Gamma)$ the  {\em medial quiver} of internal arcs of
$\Gamma$. That is, $Q(\Gamma)$ is the quiver whose set of points is the
set of internal arcs of $\Gamma$, and the arrows are defined as follows:
whenever there is a triangle $T$ in $\Gamma$ containing two internal
arcs $a$ and $b$, then $Q(\Gamma)$ contains an arrow $a \to b$ if $a$ is
a predecessor of $b$ with respect to clockwise orientation at the joint
vertex of $a$ and $b$ in $T$ (here we use that $S$ is an oriented
surface, this allows
to talk about clockwise orientation around each marked point).
\bigskip

{\bf Example.}
We illustrate the construction of $Q(\Gamma)$ when $\Gamma$ is a
triangulation of an octagon:

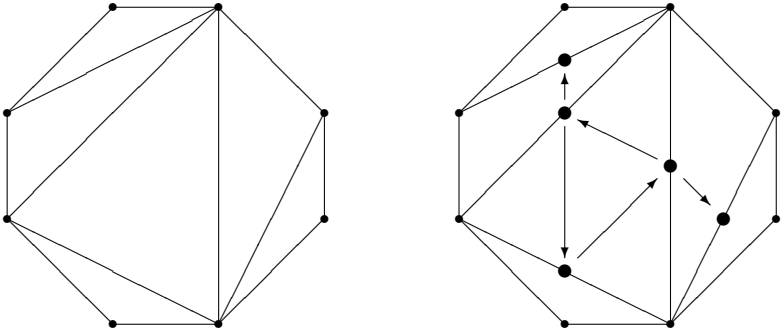
\begin{center}
\begin{picture}(120,120)(0,0)
   \put(0,40){\vertex}
   \put(0,80){\vertex}
   \put(40,0){\vertex}
  \put(40,120){\vertex}
 \put(80,0){\vertex}
  \put(80,120){\vertex}
   \put(120,40){\vertex}
   \put(120,80){\vertex}
  
    \put(40,0){\line(1,0){40}}
    \put(40,120){\line(1,0){40}}
    \put(0,40){\line(0,1){40}}
    \put(120,40){\line(0,1){40}}
    \put(40,0){\line(-1,1){40}}
   \put(120,80){\line(-1,1){40}}
   \put(0,80){\line(1,1){40}}
   \put(80,0){\line(1,1){40}}
  
   \put(0,80){\line(2,1){80}}
   \put(0,40){\line(1,1){80}}
   \put(0,40){\line(2,-1){80}}
   \put(80,0){\line(0,1){120}}
   \put(80,0){\line(1,2){40}}
   \end{picture}
   \qquad\qquad
\begin{picture}(120,120)(0,0)
   \put(0,40){\vertex}
   \put(0,80){\vertex}
   \put(40,0){\vertex}
  \put(40,120){\vertex}
 \put(80,0){\vertex}
  \put(80,120){\vertex}
   \put(120,40){\vertex}
   \put(120,80){\vertex}
   \put(40,100){\bvertex}
    \put(40,80){\bvertex}
    \put(80,60){\bvertex}
    \put(100,40){\bvertex}
    \put(40,20){\bvertex}
   \put(40,85){\vector(0,1){10}}
   \put(75,62.5){\vector(-2,1){30}}
   \put(40,75){\vector(0,-1){50}}
   \put(45,25){\vector(1,1){30}}
   \put(85,55){\vector(1,-1){10}}
  
    \put(40,0){\line(1,0){40}}
    \put(40,120){\line(1,0){40}}
    \put(0,40){\line(0,1){40}}
    \put(120,40){\line(0,1){40}}
    \put(40,0){\line(-1,1){40}}
   \put(120,80){\line(-1,1){40}}
   \put(0,80){\line(1,1){40}}
   \put(80,0){\line(1,1){40}}

   \put(0,80){\line(2,1){80}}
   \put(0,40){\line(1,1){80}}
   \put(0,40){\line(2,-1){80}}
   \put(80,0){\line(0,1){120}}
   \put(80,0){\line(1,2){40}}
  
\end{picture}
\hfill .
\end{center}
\bigskip

\begin{lemma}\label{no2cycle}
The quiver $Q(\Gamma)$ contains no oriented cycles of length $\le 2$.
\end{lemma}
{\bf Proof.}
We first show that  $Q(\Gamma)$ contains no loops.
A loop $\alpha$ at the point $a$ of $Q(\Gamma)$ would arise from a
triangle $T$ in $\Gamma$ in the following way:
$$
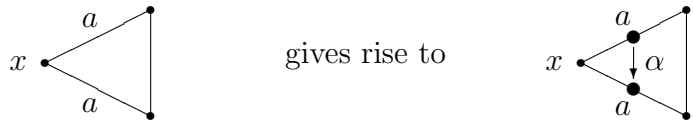
\begin{picture}(40,20)(0,20)
   \put(0,20){\vertex}
   \put(40,0){\vertex}
   \put(40,40){\vertex}
   \put(0,20){\line(2,1){40}}
   \put(0,20){\line(2,-1){40}}
   \put(40,40){\line(0,-1){40}}
    \put(-10,20){\makebox(0,0){$x$}}
    \put(17,36){\makebox(0,0){$a$}}
    \put(17,4){\makebox(0,0){$a$}}
\end{picture}
\qquad\qquad \mbox{ gives rise to } \qquad\qquad
\begin{picture}(40,20)(0,20)
   \put(0,20){\vertex}
   \put(40,0){\vertex}
   \put(40,40){\vertex}
     \put(20,30){\bvertex}
    \put(20,10){\bvertex}
   \put(0,20){\line(2,1){40}}
   \put(0,20){\line(2,-1){40}}
   \put(40,40){\line(0,-1){40}}
      \put(20,26){\vector(0,-1){12}}
    \put(-10,20){\makebox(0,0){$x$}}
    \put(28,20){\makebox(0,0){$\alpha$}}
    \put(16,37){\makebox(0,0){$a$}}
    \put(16,3){\makebox(0,0){$a$}}
\end{picture}$$
\bigskip

But in this case the triangle $T$ is homeomorphic to
$$
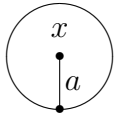
\begin{picture}(40,40)(0,0)
   \put(20,20){\vertex}
   \put(20,20){\circle{40}}
   \put(20,0){\vertex}
    \put(20,20){\line(0,-1){20}}
   \put(20,30){\makebox(0,0){$x$}}
    \put(25,10){\makebox(0,0){$a$}}
\end{picture}$$

which means that $x$ is an internal vertex, contradicting our assumption
that $M$ is contained in the boundary of $S$.

We now show that  $Q(\Gamma)$ contains no oriented cycles of length two.
Indeed, such a cycle corresponds to the following situation in $\Gamma$:
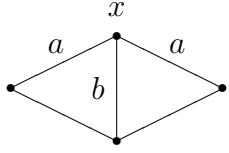
\bigskip

\begin{center}
\begin{picture}(80,40)(0,0)
   \put(0,20){\vertex}
   \put(80,20){\vertex}
   \put(40,0){\vertex}
   \put(40,40){\vertex}
   \put(0,20){\line(2,1){40}}
   \put(40,0){\line(2,1){40}}
   \put(0,20){\line(2,-1){40}}
   \put(40,40){\line(2,-1){40}}
   \put(40,40){\line(0,-1){40}}
    \put(40,50){\makebox(0,0){$x$}}
    \put(17,36){\makebox(0,0){$a$}}
    \put(63,36){\makebox(0,0){$a$}}
    \put(33,20){\makebox(0,0){$b$}}
\end{picture}
\end{center}
Then a neighbourhood of  $x$  is homeomorphic to
$$
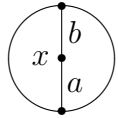
\begin{picture}(40,40)(0,0)
   \put(20,20){\vertex}
   \put(20,20){\circle{40}}
   \put(20,0){\vertex}
  \put(20,40){\vertex}
    \put(20,20){\line(0,-1){20}}
    \put(20,20){\line(0,1){20}}
   \put(12,20){\makebox(0,0){$x$}}
    \put(25,10){\makebox(0,0){$a$}}
    \put(25,30){\makebox(0,0){$b$}}
\end{picture}$$
which again contradicts the assumption that $\Gamma$ contains no
internal vertices.
\hfill $\Box$ \bigskip

{\bf Remark.} In \cite{FST} the authors associate a skew-symmetric matrix $B(\Gamma)$  to a triangulation $\Gamma$ of $(S,M)$. This construction is equivalent to the construction of  the quiver  $Q(\Gamma)$ we consider here. Since $Q(\Gamma)$ contains no oriented cycles of length $\le 2$, it is uniquely determined by a skew-symmetric matrix $B$ (where the number of arrows between two vertices is given by the entries of $B$, and the direction of the arrows is determined by the sign of the matrix entries).
It is easy to see that $B$ coincides with  $B(\Gamma)$.
Thus all the results from \cite{FST} apply, in particular,
 mutations of the quiver $Q(\Gamma)$ correspond to flips of the
triangulation $\Gamma$:

Let $b$ be an internal arc of $\Gamma$. Thus $b$ is one diagonal of the
quadrilateral formed by the  two triangles of $\Gamma$ that contain $b$.
The flip of $b$ replaces the edge $b$ by the other diagonal $b^*$ of the
same quadrilateral. Keeping all other edges unchanged, one obtains a new
triangulation $\mu_b (\Gamma)$.

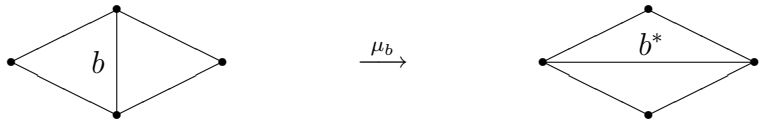
\begin{center}
\begin{picture}(80,40)(0,17)
   \put(0,20){\vertex}
   \put(80,20){\vertex}
   \put(40,0){\vertex}
   \put(40,40){\vertex}
   \put(0,20){\line(2,1){40}}
   \put(40,0){\line(2,1){40}}
   \put(0,20){\line(2,-1){40}}
   \put(40,40){\line(2,-1){40}}
   \put(40,40){\line(0,-1){40}}
    \put(33,20){\makebox(0,0){$b$}}
\end{picture}
\qquad\qquad  
$\stackrel{\mu_b}{\longrightarrow}$
\qquad\qquad
\begin{picture}(80,40)(0,17)
   \put(0,20){\vertex}
   \put(80,20){\vertex}
   \put(40,0){\vertex}
   \put(40,40){\vertex}
   \put(0,20){\line(2,1){40}}
   \put(40,0){\line(2,1){40}}
   \put(0,20){\line(2,-1){40}}
   \put(40,40){\line(2,-1){40}}
   \put(0,20){\line(1,0){80}}
    \put(41,27){\makebox(0,0){$b^*$}}
\end{picture}
\end{center}
\bigskip\bigskip

An essential ingredient in the definition of cluster algebras by
Fomin and Zelevinsky \cite{FZ}  is the mutation of skew-symmetric matrices. Reformulated in the language of quivers, one obtains a mutation of quivers $Q \mapsto \mu_b(Q)$.
The following proposition shows that  flips of the triangulation
commute with quiver mutations:

\begin{proposition}[{\cite[Prop 4.8]{FST}}]\label{flip}
Suppose that the triangulation $\mu_b(\Gamma)$ is obtained from $\Gamma$
by a flip replacing the diagonal labelled $b$.
Then $$Q(\mu_b(\Gamma)) = \mu_b(Q(\Gamma))$$
\end{proposition}


\subsection{The algebra $A(\Gamma)$}

We define in this section an algebra $A(\Gamma)$ for each triangulation
$\Gamma$ of the unpunctured marked surface $(S,M)$. Our construction
generalizes the one given in \cite{CCS} for polygons. An even more
general case is considered by Labardini in \cite{La}, where such an
algebra $A(\Gamma)$ is defined for a general marked surface (allowing
punctures). If there are no punctures,  the definitions coincide (although Labardini works in the equivalent framework of opposite medial quivers).

A triangle $T$ in $\Gamma$ is called an {\em internal triangle } if all edges
of $T$ are internal arcs. Every internal triangle $T$ in $\Gamma$ gives
rise to an oriented cycle $\alpha_T \beta_T \gamma_T$ in $Q(\Gamma)$,
unique up to cyclic permutation of the factors $\alpha_T, \beta_T, \gamma_T$.
We define\vspace{-2mm}
 $$W= \sum_{T} \alpha_T \beta_T \gamma_T$$\vspace{-2mm}
 
 where the sum runs over all internal triangles $T$ of $\Gamma$.
 Then $W$ is a potential on $Q(\Gamma)$ and we define $A(\Gamma)$ to be
the (non-completed) Jacobian algebra of $(Q,W)$ (see \cite{DWZ}, \cite{Ke}).
 Thus $A(\Gamma)$ can be described as a quotient $A(\Gamma) =
kQ(\Gamma)/I(\Gamma)$ of the path algebra $kQ(\Gamma)$ by the ideal
$I(\Gamma)$ generated by all paths $\alpha_T \beta_T$, $\beta_T\gamma_T$
and $\gamma_T\alpha_T $ whenever $T$ is an internal triangle of $\Gamma$.
In \cite[Theorem 30]{La}  it is shown that flips in the triangulation correspond to
mutations of the quiver with potential $(Q(\Gamma), W)$ as defined in
\cite{DWZ}.
\medskip

The following result is shown in \cite[Theorem 36]{La} for the more
general case of punctured marked surfaces.
\begin{lemma}\label{findim}
Let $\Gamma$ be a triangulation of  an unpunctured marked surface
$(S,M)$. Then the algebra $A(\Gamma)$ is finite-dimensional.
\end{lemma}

We show in Lemma \ref{gentle} that the algebras $A(\Gamma)$ belong to a class of algebras called {\em gentle} algebras. 
Recall from \cite{AS} that a finite-dimensional  algebra is  {\em  gentle} if it admits a presentation $A=kQ/I$ satisfying the following {\nobreak conditions:}
\begin{itemize}
\item[(G1)] At each point of $Q$ start at most two arrows and stop at
  most two arrows.\vspace{-2mm}
\item[(G2)] The ideal $I$ is generated by paths of length 2.\vspace{-2mm}
\item[(G3)] For each arrow $\beta$ there is at most one arrow $\alpha$
  and at most one arrow $\gamma$ such that $\alpha \beta \in I$
  and $\beta \gamma \in I$.\vspace{-2mm}
\item[(G4)] For each arrow $\beta$ there is at most one arrow $\alpha$
  and at most one arrow $\gamma$ such that $\alpha \beta \not\in I$
  and $\beta \gamma \not\in I$.
\end{itemize}

If the pair $(Q,I)$ satisfies conditions (G1) to (G4), we call it a
{\em gentle bound quiver}, or a {\em gentle presentation} of $A=kQ/I$.
Note that in contrast to \cite{AS}, we do not assume that $A=kQ/I$ is
triangular.
An algebra $A=kQ/I$ where $I$ is generated by paths and $(Q,I)$ satisfies  the two conditions (G1) and (G4) is called a {\em string algebra} (see \cite{BR}), thus every gentle algebra is a string algebra.
The gentle algebras can be characterized  by the fact that
their repetitive categories are special biserial
\cite{AS,PoS}.
\medskip

We recall here the classification of indecomposable modules over a string algebra $A=kQ/I$ which is given in \cite{BR} in terms   of reduced walks in the quiver $Q$:
A {\em string} in $A$ is by definition a reduced walk $w$ in $Q$ avoiding the zero-relations, thus $w$ is a sequence
$$ w = \;x_1 \stackrel{\alpha_1}{\longleftrightarrow} x_2 \stackrel{\alpha_1}{\longleftrightarrow}\cdots  \stackrel{\alpha_{n-1}}{\longleftrightarrow} x_n
$$
where the $x_i$ are vertices of $Q$ and  each $\alpha_i$ is an arrow between the vertices $x_i$ and $x_{i+1}$ in either direction such that $w$ does not contain a sequence of the form 
$ \stackrel{\beta}{\longleftarrow} \stackrel{\beta}{\longrightarrow}$ or $ \stackrel{\beta_1}{\longrightarrow} \cdots\stackrel{\beta_s}{\longrightarrow}$ with $\beta_1 \cdots \beta_s \in  I$ or their duals. A string is {\em cyclic} if the first and the last vertex coincide.
A {\em band} is defined to be a cyclic string $b$ such that each power $b^n$ is a string, but $b$ itself is not a proper power of some string $c$.
\bigskip

The {\em string module} $M(w)$ is obtained from the string $w$ by replacing each $x_i$ in $w$ by a copy of the field $k$. The action of an arrow $\alpha$ on $M(w)$ is induced by the relevant identity morphisms if
$\alpha$ lies on $w$, and is zero otherwise.
The dimension vector $\udim M(w)$ of $M(w)$ is obtained by counting how often the string $w$ passes through each vertex $x$ of the quiver $Q$:
$$ \udim M(w) = ( \sum_{1 \le i \le n} \delta_{x,x_i}  )_{x \in Q_0}
$$
where $ \delta_{x,x_i} =1$ for $x=x_i$ and  $ \delta_{x,x_i} =0$ otherwise.
Similarly, each band $b$ in $A$ gives rise to a family of {\em band modules} $M(b,\lambda, n)$ where $\lambda \in k$ and $n \in \mN$ (we refer to \cite{BR} for the precise definition).
All string and band modules are indecomposable, and in fact every indecomposable $A-$module is either a string module $M(w)$  or a band module $M(b,\lambda, n)$, see  \cite{BR}.
\newpage

We now return to the study of algebras stemming from surface triangulations:
\begin{lemma}\label{gentle}
Let $\Gamma$ be a triangulation of  an unpunctured marked surface
$(S,M)$. Then $A(\Gamma)$ is a gentle algebra.
\end{lemma}
{\bf Proof.}
By Lemma \ref{findim}, the algebra $A(\Gamma)$ is finite-dimensional, so
we only need to verify conditions (G1) to (G4) for the bound quiver
$(Q(\Gamma), I(\Gamma))$  of $A(\Gamma)$.

(G2) By definition, the ideal $I(\Gamma)$ is generated by paths of
length two.
\smallskip

(G1) Let $a$ be a point of $Q(\Gamma)$ corresponding to an internal arc $a$ of
$\Gamma$.
Since $\Gamma$ is a triangulation of a surface, the arc $a$ is
contained in at most two triangles:

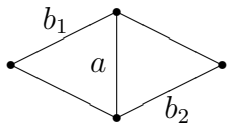
\begin{center}
\begin{picture}(80,40)(0,0)
   \put(0,20){\vertex}
   \put(80,20){\vertex}
   \put(40,0){\vertex}
   \put(40,40){\vertex}
   \put(0,20){\line(2,1){40}}
   \put(40,0){\line(2,1){40}}
   \put(0,20){\line(2,-1){40}}
   \put(40,40){\line(2,-1){40}}
   \put(40,40){\line(0,-1){40}}
    \put(17,37){\makebox(0,0){$b_1$}}
    \put(63,3){\makebox(0,0){$b_2$}}
    \put(33,20){\makebox(0,0){$a$}}
\end{picture}
\end{center}
Hence there are at most two arrows $\alpha_1: b_1 \to a$ and $\alpha_2:
b_2 \to a$ of $Q(\Gamma)$ ending in $a$. The same holds for arrows
starting in a point $a$.
\medskip

(G3),(G4)  Suppose now that $Q(\Gamma)$ contains $\alpha_1, \alpha_2,
\beta$ as follows:
\bigskip

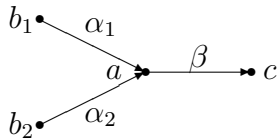
\begin{center}
\begin{picture}(80,40)(0,0)
   \put(40,20){\vertex}
   \put(80,20){\vertex}
   \put(0,0){\vertex}
   \put(0,40){\vertex}
    \put(0,0){\vector(2,1){40}}
   \put(0,40){\vector(2,-1){40}}
   \put(40,20){\vector(1,0){40}}
    \put(23,37){\makebox(0,0){$\alpha_1$}}
    \put(23,3){\makebox(0,0){$\alpha_2$}}
    \put(60,25){\makebox(0,0){$\beta$}}
   \put(-7,40){\makebox(0,0){$b_1$}}
   \put(-7,0){\makebox(0,0){$b_2$}}
   \put(28,20){\makebox(0,0){$a$}}
  \put(87,20){\makebox(0,0){$c$}}
 \end{picture}
\end{center}

We have to show that precisely one of $\alpha_1\beta$,$\alpha_2\beta$
belongs to $I(\Gamma)$.
In $\Gamma$, the internal arcs $a,b_1,b_2$  belong to two triangles as
considered in the proof of (G1). The arrow $\beta$ encodes that the
arc $c$  is a successor of $a$ in one of these triangles, say the one
formed by  $a,b_1,c$. This gives rise to the relation $\alpha_1\beta$,
and $\alpha_2\beta$ does not belong to $I(\Gamma)$ since $\alpha_2$ and
$\beta$ arise from different triangles.
\hfill $\Box$ \bigskip

From the construction of $A(\Gamma)$ it is clear that for each
$\alpha\beta \in I(\Gamma)$ there is an arrow $\gamma$ in $Q(\Gamma)$
such that $\beta\gamma \in I(\Gamma)$ and $\gamma \alpha \in I(\Gamma)$.
In the following lemma we study a homological property of all gentle
algebras satisfying this condition: An algebra $A$ is {\em Gorenstein of dimension one} if the injective
dimension of the (finitely generated) projective $A-$modules is at most one, and the
projective dimension of the (finitely generated) injective $A-$modules is at most one.
Note that all cluster-tilted algebras are Gorenstein of dimension one (see 3.1 below), and that an algebra of Gorenstein dimension one is either hereditary or has infinite global dimension (as shown in \cite{KR}).

\begin{lemma}\label{gorenstein}
Let $A=kQ/I$ be a gentle algebra such that for each $\alpha\beta \in I$
there is an arrow $\gamma$ in $Q$ such that $\beta\gamma \in I$ and
$\gamma \alpha \in I$.
Then $A$ is Gorenstein of dimension one.
\end{lemma}
{\bf Proof.}
We only compute the projective
dimension of the injective modules here, the proof of the other part in the definition of Gorenstein of dimension one is dual. 
It is sufficient to show that for every vertex $x$ of $Q$ the corresponding
indecomposable injective $A-$module $I_x$ has projective dimension at
most one. To do so, we construct explicitly a projective resolution of
$I_x$. We write the string module $I_x$ as $I_x=
M(u_1\alpha_1\alpha_2^{-1}u_2^{-1})$ where $u_1$ and $u_2$ are oriented
paths.  Both paths might have length zero, and in this case, also the
arrows $\alpha_1$ and $\alpha_2$ might not be present.
The following figure is used throughout the proof:
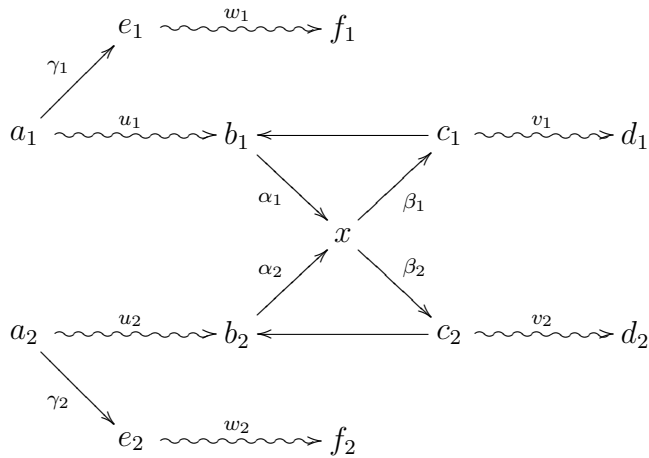
\begin{displaymath}
        \xymatrix{
                   & e_1\ar@{~>}[rr]^{w_1} & & f_1 & &        \\
                   a_1\ar@{~>}[rr]^{u_1}\ar[ur]^{\gamma_1} & &
b_1\ar[dr]_{\alpha_1} & & c_1\ar[ll]\ar@{~>}[rr]^{v_1} & & d_1\\
                    & & & x\ar[ur]_{\beta_1}\ar[dr]^{\beta_2} & & & \\
                   a_2\ar@{~>}[rr]^{u_2}\ar[dr]_{\gamma_2} & &
b_2\ar[ur]^{\alpha_2} & & c_2\ar[ll]\ar@{~>}[rr]^{v_2} & & d_2\\
  & e_2\ar@{~>}[rr]^{w_2} & & f_2 & &        \\        }
\end{displaymath}
Note that $\{x,c_1,b_1\}$ and $\{x,c_2,b_2\}$ form oriented cycles in
$Q$ such that the composition of any two consecutive arrows is zero.
Let $p_0:P(0) \to I_x$ be a projective cover, then
$$ P(0) = M(w_1^{-1}\gamma_1^{-1}u_1\alpha_1\beta_2v_2) \oplus
M(w_2^{-1}\gamma_2^{-1}u_2\alpha_2\beta_1v_1)$$
and
$$\Ker p_0 = M(w_1) \oplus M(w_2)\oplus M(v_1^{-1}\beta_1^{-1}\beta_2v_2)$$
(note that some summands of the terms of this sequence can be zero).
We show that $\Ker p_0$ is projective, thus obtaining the desired
projective resolution
\begin{displaymath}
        \xymatrix{ 0 \ar[r] & \Ker p_0 \ar[r] &  P(0)\ar[r]^{p_0} &
I_x\ar[r] & 0
        }
\end{displaymath}
In order to see that the first two summands of $\Ker p_0$ are projective
(namely the indecomposable projectives $P_{e_1}$ and $P_{e_2}$) one has to show that there are no other
arrows starting at the vertices $e_1,e_2$.
Suppose there is an arrow $\delta_1 : e_1 \to y$ in $Q$. Since the
algebra $A$ is gentle, the composition $\gamma_1\delta_1$ lies in the
ideal $I$. The assumption of the lemma guarantees the existence of a
cycle  $\gamma_1\delta_1\epsilon_1$  such that
$\gamma_1\delta_1,\delta_1\epsilon_1,\epsilon_1\gamma_1 \in I$. 
But then the simple $A-$module $S_y$ would be a composition factor of $I_x$, contradicting the
 assumption $I_x= M(u_1\alpha_1\alpha_2^{-1}u_2^{-1})$.
 This shows that $M(w_1)=P_{e_1}$, and a similar argument shows that $M(w_2)=P_{e_2}$.
 Since $M(v_1^{-1}\beta_1^{-1}\beta_2v_2)=P_x$, we conclude that $\Ker p_0$ is projective.
\hfill $\Box$
\bigskip

{\bf Example.}
We illustrate the projective resolution constructed in Lemma \ref{gorenstein} when $\Gamma$ is the following triangulation of a polygon with 11 vertices (where the midpoints of internal arcs are labeled):
\bigskip

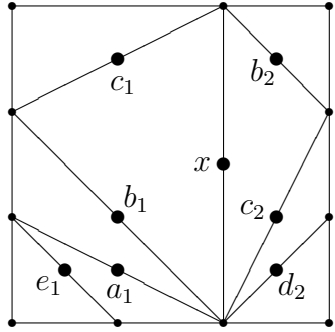
\begin{center}
\begin{picture}(120,120)(0,0)
   \put(0,0){\vertex}
   \put(40,0){\vertex}
   \put(80,0){\vertex}
  \put(120,0){\vertex}
 \put(0,40){\vertex}
  \put(120,40){\vertex}
   \put(0,80){\vertex}
   \put(120,80){\vertex}
   \put(0,120){\vertex}
    \put(80,120){\vertex}
    \put(120,120){\vertex}
    \put(100,40){\bvertex}
     \put(100,100){\bvertex}
    \put(100,20){\bvertex}
    \put(20,20){\bvertex}
   \put(80,60){\bvertex}
    \put(40,20){\bvertex}
     \put(40,40){\bvertex}
    \put(40,100){\bvertex}
      
    \put(0,0){\line(1,0){120}}
    \put(0,0){\line(0,1){120}}
    \put(120,0){\line(0,1){120}}
    \put(0,120){\line(1,0){120}}
    \put(40,0){\line(-1,1){40}}
   \put(120,80){\line(-1,1){40}}
   \put(80,0){\line(1,1){40}}

   \put(0,80){\line(2,1){80}}
   \put(80,0){\line(-1,1){80}}
   \put(80,0){\line(-2,1){80}}
   \put(80,0){\line(0,1){120}}
   \put(80,0){\line(1,2){40}}
  
   \put(72,60){\makebox(0,0){$x$}}
   \put(95,95){\makebox(0,0){$b_2$}}
    \put(91,43){\makebox(0,0){$c_2$}}
   \put(106,14){\makebox(0,0){$d_2$}}
  \put(14,14){\makebox(0,0){$e_1$}}
  \put(41,11){\makebox(0,0){$a_1$}}
 \put(47,47){\makebox(0,0){$b_1$}}
 \put(42,90){\makebox(0,0){$c_1$}}
 
\end{picture}
\end{center}
\bigskip

The corresponding algebra $A(\Gamma)$ is given by the quiver

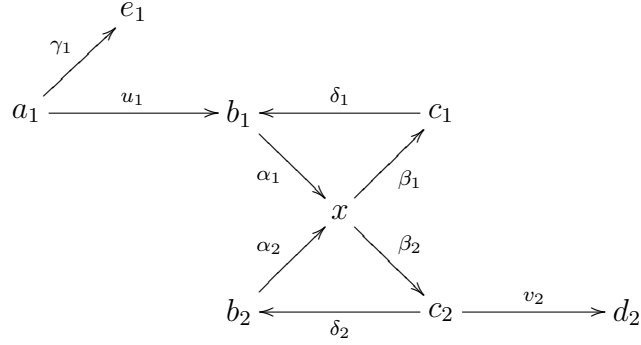
\begin{displaymath}
        \xymatrix{
                   & e_1  & &  & &        \\
                   a_1\ar[rr]^{u_1}\ar[ur]^{\gamma_1} & &
b_1\ar[dr]_{\alpha_1} & & c_1\ar[ll]_{\delta_1} & & \\
                    & & & x\ar[ur]_{\beta_1}\ar[dr]^{\beta_2} & & & \\
                    & & b_2\ar[ur]^{\alpha_2} & & c_2\ar[ll]^{\delta_2}\ar[rr]^{v_2} & & d_2\\    }
\end{displaymath}

with relations $\alpha_1\beta_1= \beta_1\delta_1=\delta_1\alpha_1 = 0$ and $\alpha_2\beta_2= \beta_2\delta_2=\delta_2\alpha_2 = 0$.
The projective resolution of the injective module $I_x$ is then
\begin{displaymath}
        \xymatrix{ 0 \ar[r] & P_{e_1} \oplus P_x \ar[r] &  P_{a_1} \oplus P_{b_2} \ar[r] &
I_x\ar[r] & 0
        }
\end{displaymath}
where $P_{e_1}$ is simple and $P_x = M(\beta_1^{-1}\beta_2v_2)$, 
$P_{a_1}= M(u_1\alpha_1\beta_2v_2) $ and $P_{b_2}=M(\alpha_2\beta_1)$.
\bigskip

We recall from \cite{Ga} the concept of Galois coverings of  bound quiver algebras:
Let $\Lambda=k\tilde{Q}/\tilde{I}$ be a bound quiver algebra (where the quiver $\tilde{Q}$ is not necessarily finite). A group $G$ of $k-$linear automorphisms of $\Lambda$ is acting freely on $\Lambda$ if $ge_x \neq e_x$ for each vertex $x$ of $\tilde{Q}$ and each $g \neq 1 $ in $G$. 
In this case the multiplication in $\Lambda$ induces a multiplication on the set $\Lambda/G$ of $G-$orbits which turns $\Lambda/G$ into an algebra. The canonical projection $\Lambda \to \Lambda/G$ is called the {\em Galois covering of $\Lambda/G$ with group $G$.}

In the following theorem we call (as in \cite{BV})  {\em 3-cycle}  an oriented cycle $\alpha\beta\gamma$ where $\alpha, \beta, \gamma$ are three distinct arrows, and by a {\em 3-cycle with radical square zero} we mean a
3-cycle $\alpha\beta\gamma$ in an algebra $kQ/I$ such that $\alpha\beta,
\beta\gamma,\gamma\alpha \in I$.
Moreover, by a {\em simple cycle} we refer to a subquiver $C$ of $Q$ with $n$ distinct vertices $\{ x_0, x_1, \ldots , x_{n-1}, x_n=x_0 \}$ and $n$ arrows $\alpha_i: x_i \to x_{i+1}$, for $i= 1, \ldots, n-1$. 

\begin{theorem}\label{structure}
Let $\Gamma$ be a triangulation of an unpunctured marked  surface
$(S,M)$.  Then:

\begin{itemize}
\item[(1)] the algebra $A(\Gamma)$ is a gentle algebra,

\item[(2)] the algebra $A(\Gamma)$ is Gorenstein of dimension one,

\item[(3)]  there is a relation in $A(\Gamma)$ from $x$ to $y$ only if there is 
an arrow $y \to x$,

\item[(4)] $A(\Gamma)$ admits a Galois covering by a bound quiver algebra $k\tilde{Q}/\tilde{I}$ satisfying:
\begin{itemize}
\item[(T1)] Every simple cycle in $\tilde{Q}$ is a 3-cycle with radical
square zero,
\item[(T2)] The only relations in $\tilde{I}$ are those in the 3-cycles.
\end{itemize}\end{itemize}

\end{theorem}
{\bf Proof.}
Part (1) is shown in Lemma \ref{gentle}, and Part (2) is shown in Lemma \ref{gorenstein} since the condition imposed
on the gentle algebra $A$ there clearly holds for the algebra $A(\Gamma)$.
Part (3)  follows directly from the definition of $A(\Gamma)$.
Maybe the most intuitive way to obtain the Galois covering required in Part (4) is the following: 
By construction, the only relations in the algebra $A(\Gamma)$ are those in the 3-cycles. In a first step, we identify all 3-cycles to points, replacing each 3-cycle $C$ with vertices $\{x_1,x_2,x_3\}$ by one single vertex $x$ and replacing each arrow $y \to x_i$ (or $x_i \to y$, respectively) by an arrow $y \to x$ (or $x \to y$, respectively).
The  quiver $\overline{Q}$ thus obtained contains no relations, and we let $\widetilde{\overline{Q}}$ be its universal Galois covering, a (maybe infinite) tree.
The bound quiver $(\tilde{Q},\tilde{I})$ is then obtained by placing back the 3-cycles $C=\{x_1,x_2,x_3\}$ for all contracted vertices $x$ of $\widetilde{\overline{Q}}$.
\hfill $\Box$ \bigskip

Note that the finite quivers satisfying conditions (T1) and (T2) from the previous theorem form precisely the
class of quivers ${\cal Q}_n$ considered in \cite{BV}, where also  the same relations
are imposed. It would be interesting to relate the Galois covering $(\tilde{Q},\tilde{I})$ constructed above with the universal cover of the bordered surface $(S,M)$.


\subsection{Recovering topological data from $A(\Gamma)$}
The condition (4) in Theorem \ref{structure} is very strong. Combined
with the fact that the algebra is gentle, it implies the
remaining conditions (2) and (3).
We show in this section that a gentle algebra  satisfying
condition (4) is given by  an unpunctured marked surface.
 \bigskip
 
First we give a different combinatorial description of the algebras
studied here. Consider the following two bound quivers, where type I is
a quiver of type $\mA_2$, and type II is a 3-cycle with radical square zero:

$$\mbox{ Type I }\qquad  
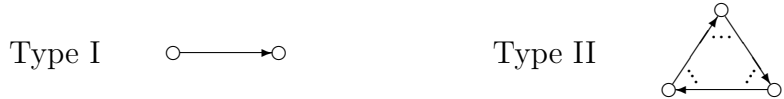
\begin{picture}(70,0)(0,-4)
   \put(0,0){\bwertex}
   \put(40,0){\bwertex}
    \put(2.5,0){\vector(1,0){35}}
  \end{picture}
\qquad\qquad \mbox{ Type II } \qquad
\begin{picture}(40,30)(0,10)
   \put(0,0){\bwertex}
   \put(40,0){\bwertex}
   \put(20,30){\bwertex}
\put(23,20){\svertex}
\put(20,20){\svertex}
\put(17,20){\svertex}
\put(29.66,3){\svertex}
\put(31,5){\svertex}
\put(32.33,7){\svertex}
\put(10.33,3){\svertex}
\put(9,5){\svertex}
\put(7.66,7){\svertex}
\put(1.5,2.25){\vector(2,3){17}}
\put(21.5,27.75){\vector(2,-3){17}}
\put(37.5,0){\vector(-1,0){35}}
\end{picture}$$
\bigskip

Using these bound quivers one can construct algebras in the following
way: Suppose we start with a collection $C$ of disjoint blocks of type I
and II. Choose a partial matching (that is to say a partial bijection) $\pi$ of the vertices in $C$, where
matching a vertex to  another vertex of the same block is
not allowed. Identifying (or ''gluing'') the vertices within each pair
of the matching we obtain an algebra $A(C,\pi)$. Note that the arrows are not identified by this procedure, so one might obtain parallel arrows or two-cycles. 
We consider only
matchings where the algebra $A(C,\pi)$ is connected. 

\medskip

The procedure of gluing blocks is considered  in a more general situation (using plenty of building blocks) in \cite{Br}, where the
resulting  algebras  are called {\em kit algebras}. A similar construction to glue blocks of type
I, II and four more types is  described in \cite{FST}.
\medskip

We show below that  the gentle algebras that admit a Galois covering  satisfying conditions (T1) and (T2) from Theorem \ref{structure} are algebras of the form $A(C,\pi)$, thus results from \cite{FST} concerning these algebras can be applied.

\begin{proposition}\label{topology}
Let $A = kQ/I$ be a gentle algebra where every relation lies on a
3-cycle with radical square zero.  Then there exists an unpunctured marked
surface $(S,M)$ with a triangulation $\Gamma$ such that $A(\Gamma) = A$.
\end{proposition}

{\bf Proof.} The statement follows from \cite[(14.1)]{FST} once we show that the algebra $A$ admits  a unique block decomposition $A=A(C,\pi)$ using blocks of type I and II. 
We define thus $C$ to be the disjoint union of all 3-cycles with radical square zero of $A$ together with the disjoint union of all remaining arrows from $A$.
Denote by $f$ the quiver morphism $f:C \to Q$ that identifies the blocks of $C$ with their images in $Q$.
\medskip

We first show that $|f^{-1}(x)| \le 2$ for each vertex $x \in Q$.
Indeed, if $f^{-1}(x)$ contains three different vertices, then there are three different arrows in $Q$ adjacent to the vertex $x$.
But since the algebra $A$ is gentle, there has to be one relation between these three arrows. 
However, the set $C$ is constructed in such a way that all relations of $A$ belong to one of the components in $C$, so there are no relations between arrows corresponding to different components of $C$, and so the fiber $f^{-1}(x)$
 contains at most two vertices.
 
 We now define a matching $\pi$ on $C$ relating $x_1$ to $x_2$ whenever $f^{-1}(x) = \{ x_1, x_2 \}$.
 As required in the definition of $A(C,\pi)$, we do not match a vertex to itself or to some vertex in the some block. It is clear from the construction that $A= A(C,\pi)$.
 Moreover, the choice of blocks of type I or II is unique since all relations have to correspond to a block of type II.
\hfill $\Box$ \bigskip

We would like to point out that all algebras $A(\Gamma)$ given by a triangulation $\Gamma$ of an unpunctured marked surface are of the form $A(C,\pi)$ for some $C$ and $\pi$, but the converse is not true: One can easily produce two-cycles in an algebra $A(C,\pi)$, but this never occurs for the algebras $A(\Gamma)$ as we have shown in Lemma \ref{no2cycle}.

%
%

\section{Gentle cluster-tilted algebras}
\subsection{Cluster-tilted algebras}
Let $\Delta$ be an acyclic quiver. In \cite{BMRRT} the cluster category
${\cal C}_{\Delta}$ is studied in order to obtain a categorical
interpretation of the cluster variables of the cluster algebra
associated with $\Delta$.
It is shown in \cite{BMRRT} that clusters correspond bijectively to
tilting objects $T$ in ${\cal C}_{\Delta}$.
Their endomorphism rings $\End_ {{\cal C}_\Delta}(T)$ are called {\em cluster-tilted algebras of type $\Delta$. }They  were introduced and studied in
\cite{BMR}.
\medskip

We use here a  different description that has been  given in \cite{ABS}.
Denote by $A$ the hereditary algebra $A=k\Delta$.
An $A$-module $T$ is called a {\em tilting module} provided
$\Ext{1}{A}(T,T)=0$ and the number of
isomorphism classes of indecomposable summands of $ T$ equals
the number of isomorphism classes of simple $A$-modules.
In this case, the endomorphism ring $\End_A(T)$ is called a {\em tilted
algebra of type $\Delta$.}
\medskip

Let $C$ be an algebra of global dimension two.  The trivial extension
\begin{displaymath}
        \tilde{C} = C \ltimes \Ext{2}{C}(DC, C)
\end{displaymath}
of $C$ by the $C-C-$bimodule  $\Ext{2}{C}(DC, C)$ is called the
\emph{relation-extension} of $C$.
It is useful to describe explicitly the operations on $\tilde{C}$.
As an abelian group,  $  \tilde{C} = C \oplus \Ext{2}{C}(DC, C)$.
Let thus $(c,{\bf e})$ and $(c',{\bf e'})$ be two elements of  $  \tilde{C} $, where $\bf e$ and $\bf e'$ are respectively 
represented by the exact sequences of $C-$modules
\begin{displaymath}
 \xymatrix{ {\bf e}: & 0 \ar[r] & P \ar[r] &
M \ar[r] & N\ar[r] & I\ar[r] & 0\\
{\bf e'}: & 0 \ar[r] & P' \ar[r] &
M' \ar[r] & N'\ar[r] & I' \ar[r] & 0 
  }         
\end{displaymath}
\bigskip

with $P,P'$ projective and $I,I'$ injective.
The addition is given by 

$$ (c,{\bf e})+(c',{\bf e'}) =(c+c',{\bf e}+{\bf e'})
$$
\medskip

where the sum $c+c'$ is the ordinary sum inside $C$, while ${\bf e}+{\bf e'}$ is the Baer sum in $\Ext{2}{C}(DC, C)$ (for which we refer to any textbook of homological algebra). 
The product in $\tilde{C}$ is given by the formula
$$ (c,{\bf e})(c',{\bf e'}) =(cc',c{\bf e}'+{\bf e}c')
$$
where the product $cc'$ is the ordinary product inside $C$, while $c{\bf e}'$ and ${\bf e'}c$ are defined as follows. 
Viewing $c \in C$ as an element of $\End C_C \cong C$, then ${\bf e}_1 = c {\bf e}'$ is represented by the sequence obtained by pulling down the sequence $\bf e'$
\begin{displaymath}
 \xymatrix{ {\bf e'}: & 0 \ar[r] & P' \ar[r]\ar[d]^c &
M' \ar[r]\ar[d] & N' \ar[r] \ar@{=}[d] & I' \ar[r]\ar@{=}[d]  & 0\\
{\bf e}_1: & 0 \ar[r] & P_1 \ar[r] &
M_1 \ar[r] & N' \ar[r] & I' \ar[r] & 0 
  }        
\end{displaymath}
Similarly, viewing $c' \in C$ as an element of $\End DC_C \cong C$, then ${\bf e}_2 = {\bf e}c'$ is represented by the sequence obtained by lifting the sequence $\bf e$
\begin{displaymath}
 \xymatrix{ {\bf e}_2: & 0 \ar[r] & P \ar[r]\ar@{=}[d] &
M \ar[r]\ar@{=}[d] & N_2 \ar[r] \ar[d] & I_2 \ar[r]\ar[d]^{c'}  & 0\\
{\bf e}: & 0 \ar[r] & P \ar[r] &
M \ar[r] & N \ar[r] & I \ar[r] & 0 
  }        
\end{displaymath}

The following theorem allows to view cluster-tilted algebras as
relation-extensions of tilted algebras.
\begin{theorem}[\cite{ABS}]\label{theo::abs}
An algebra $\Lambda$ is cluster-tilted of type $\Delta$ if and only if
there exists a tilted algebra $C$ of type $\Delta$ such that $\Lambda$
is isomorphic to the relation-extension $\tilde{C}$ of $C$.
\end{theorem}

Every cluster-tilted algebra satisfies conditions (2) and (3) from
Theorem \ref{structure}, see \cite{KR}.
The bound quivers of cluster-tilted algebras of type $\mA$ are
explicitly described in \cite{BV}, Prop. 3.1. In fact they were already
described  in \cite{CCS} as the algebras $A(\Gamma)$ arising from a
triangulation of an unpunctured polygon.
The following proposition is contained in \cite{BV}, but it also follows from theorem 3.3 below.

\begin{proposition}[{\cite[(3.1)]{BV}}]
An algebra $A$ is cluster-tilted of type $\mA$ precisely when $A$ is gentle and there is
a presentation $A = kQ/I$ which satisfies conditions (T1) and (T2) from Theorem \ref{structure}.
\end{proposition}

In particular,  the cluster-tilted algebras of type $\mA$ are gentle.
We describe in the following theorem which of the gentle algebras are
cluster-tilted:

\begin{theorem}\label{theo::principal}
Let $C=kQ_C/I_C$ be a tilted algebra,  and $\tilde{C}$ be its
relation-extension.  The following are equivalent.

\begin{enumerate}
\item C is gentle;
\item C is tilted of type $\mathbb{A}$ or $\tilde{\mathbb{A}}$;
\item $\tilde{C}$ is gentle; and
\item $\tilde{C}$  is cluster-tilted of type $\mathbb{A}$ or
$\tilde{\mathbb{A}}$.
\end{enumerate}

\end{theorem}

The rest of this section is devoted to prove Theorem \ref{theo::principal}.
%
%
\subsection{Preliminaries}
        \label{sect::premieres implications}

One part of the proof of the main theorem follows from a result of
\cite{S} which says that the class of gentle algebras is stable under
tilting.

\begin{lemma}\label{lemm::1}
If a tilted algebra is gentle, then it is tilted of type $\mathbb{A}$
or $\tilde{\mathbb{A}}$.
\end{lemma}
{\bf Proof.} Let $\Delta$ be a quiver such that $C$ is tilted of type
$\Delta$.  Then there exists a tilting $C$-module  $T$ such that
$\textrm{End } T = k\Delta$.  According to \cite{S}, $k\Delta$ is a
gentle algebra.  This implies that the quiver $\Delta$ is  of type
 $\mathbb{A}$ or $\tilde{\mathbb{A}}$.
\hfill$\Box$


\begin{lemma}\label{lemm::2}
If $\tilde{C}$ is gentle, then so is $C$.
\end{lemma}
{\bf Proof.} This follows from the fact that $\tilde{C}$ is a split extension of $C$ and from \cite[(2.7)]{ACT}.\hfill$\Box$
\bigskip

\begin{lemma}\label{lemm::3}
The algebra $C$ is tilted of type $\mathbb{A}$ or $\tilde{\mathbb{A}}$ if and
only if $\tilde{C}$ is cluster-tilted of type $\mathbb{A}$ or
$\tilde{\mathbb{A}}$.
\end{lemma}

{\bf Proof.} It is clear that if $C$ is tilted of type $\mathbb{A}$ or
$\tilde{\mathbb{A}}$, then $\tilde{C}$ is cluster-tilted of type $\mathbb{A}$ or
$\tilde{\mathbb{A}}$.
Conversely, suppose that $\tilde{C}$ is cluster-tilted of type $\mathbb{A}$ or
$\tilde{\mathbb{A}}$. 
Then, by \cite{ABS} there exists a local slice $\Sigma'$ in $ \mod \tilde{C}$ such that $C'= \tilde{C}/\Ann \Sigma'$ is tilted of type $\mathbb{A}$ or $\tilde{\mathbb{A}}$. On the other hand, since $\tilde{C}= C \ltimes \Ext{2}{C}(DC, C)$, then there exists a local slice $\Sigma$ in $\mod \tilde{C}$ such that $C= \tilde{C}/\Ann \Sigma$.
Since both $\Sigma$ and $\Sigma'$ have the same underlying graph, then $C$ and $C'$ have the same type, so $C$ is tilted of type $\mathbb{A}$ or
$\tilde{\mathbb{A}}$.
\hfill$\Box$
\bigskip

The main part of the proof is concerned about the problem of showing that  $\tilde{C}$ is
gentle if $C$ is tilted of type
$\mathbb{A}$ or $\tilde{\mathbb{A}}$.
This will be done in the next subsection.

%
%
\subsection{Relation-extensions of  tilted algebras of types
$\mathbb{A}$ and $\tilde{\mathbb{A}}$}
  \label{sect::quatrieme implication}
  
Suppose that $C = kQ_C/I_C$ is tilted of type $\mathbb{A}$ or
$\tilde{\mathbb{A}}$.  
In particular, $C$ is gentle because of \cite{A} and \cite{AS}.
 Moreover, the quiver of $\tilde{C}$ is known, as are some of its
relations, namely those already in $C$ (see \cite{ABS} and \cite{ACT}).
 The aim here is to study the remaining relations of $\tilde{C}$.

First, the bound quiver of a tilted algebra of type $\mathbb{A}$ has been described in \cite{A}
and that of a tilted algebra of type $\tilde{\mathbb{A}}$ in \cite{R}.
The criterion given here is derived from \cite{HL}.
\bigskip

We recall that a {\em double-zero} in a gentle algebra is a reduced walk of the
form $\alpha\beta\omega\gamma\delta$, where $\alpha, \beta, \gamma$ and
$\delta$ are arrows such that $\alpha\beta$ and $\gamma\delta$ are
relations, while $\omega$ is a non-zero reduced walk (that is, a walk which does
not contain any relation).  Note that $\omega$ may be trivial and that in this case $\beta$ and $\gamma$ may coincide.
\bigskip

\begin{example}
{\rm The algebra
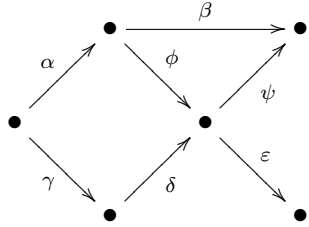
\begin{displaymath}
        \xymatrix{ & \bullet\ar[rr]^{\beta}\ar[dr]^{\phi} & & \bullet \\
                   \bullet\ar[ur]^{\alpha}\ar[dr]_{\gamma} & &
\bullet\ar[ur]_{\psi}\ar[dr]^{\varepsilon} & \\
                   & \bullet\ar[ur]_{\delta} & & \bullet
        }
\end{displaymath}
where $\alpha\beta = \phi\psi = \delta\varepsilon = 0$, is gentle with a
double-zero (namely $\phi\psi\beta^{-1}\phi\psi$).}
\end{example}
\bigskip

\begin{proposition}[{\cite{A,AS}}]\label{prop::bound quiver} 
\makebox[0mm]{}
\begin{enumerate}
  \item An algebra is tilted of type $\mathbb{A}$ if and only if it admits a bound quiver presentation $kQ/I$ with $(Q,I)$ a gentle tree with no double-zero.
  \item An algebra is tilted of type $\tilde{\mathbb{A}}$ if and only if it admits a bound quiver presentation $kQ/I$ with $(Q,I)$ a gentle presentation with no double-zero and a unique (non-oriented) cycle such that, if the cycle is a band,  then all arrows
attached to the cycle either enter it or leave it.
\end{enumerate}
\end{proposition}

\begin{example}\label{exem::an et antilde}
\rm Consider the algebras given by the bound quivers

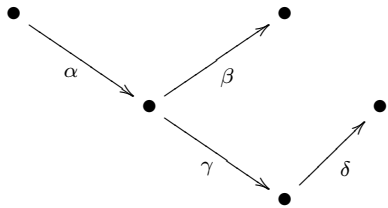
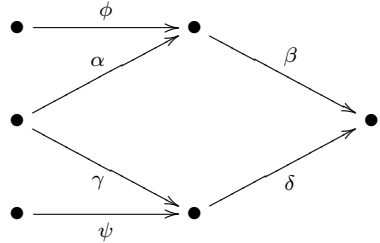
\begin{displaymath}
        \xymatrix{ \bullet\ar[dr]_{\alpha} & & \bullet &                
  & \bullet\ar[r]^{\phi} & \bullet\ar[dr]^{\beta} &   \\
                   & \bullet\ar[ur]_{\beta}\ar[dr]_{\gamma} & & \bullet
   & \bullet\ar[ur]^{\alpha}\ar[dr]_{\gamma} & & \bullet \\
                   & & \bullet\ar[ur]_{\delta} &                        
  & \bullet\ar[r]_{\psi} & \bullet\ar[ur]_{\delta} &  \\
                   & \alpha\beta = 0 & &                                
  & & \alpha\beta = \gamma\delta = 0. &
        }
\end{displaymath}
Using proposition \ref{prop::bound quiver}, we see that the first one is
tilted of type $\mathbb{A}$, while the second one is tilted of type
$\tilde{\mathbb{A}}$.
\end{example}

\subsubsection{A vanishing criterion}
We need a criterion allowing to verify whether a given exact sequence
represents the zero element in the second extension group.  We
recall the following result from \cite[(II.1.3)]{HRS}.

\begin{lemma}\label{lemm::HRS}
Given a morphism $f:M\longrightarrow N$, the exact sequence
\begin{displaymath}
        \xymatrix{0\ar[r] & \Ker f\ar[r] & M\ar[r]^{f} & N\ar[r] &
\Coker f\ar[r] & 0
        }        
\end{displaymath}
represents the zero element of $\Ext{2}{}(\Coker f, \Ker f)$ if and only
if there exist a module $X$ and morphisms $g,h$ such that the sequence
\begin{displaymath}
        \xymatrix{0\ar[r] & M\ar[r]^{(p, g)^{t}\phantom{xxxx}} & \Ima f
\oplus X \ar[r]^{\phantom{xxxx}(j,h)} & N\ar[r] & 0
        }                
\end{displaymath}
is exact, where $p$ and $j$ are the natural morphisms arising from $f$.
\end{lemma}

The following lemma will be used frequently.

\begin{lemma}\label{lemm::zero}
Let $(Q,I)$ be a gentle presentation of an algebra $C$, and let
$\alpha:c\longrightarrow b$ and $\beta:b\longrightarrow a$ be arrows in
$Q$.  Let $\sigma$ and $\eta$ be strings, not passing through $b$, such
that $\beta\sigma$ and $\eta\alpha$ are strings.  Let
$f:M(\beta\sigma)\longrightarrow M(\eta\alpha)$ be a morphism such that
$\Ima  f = S_b$.

Then the exact sequence
\begin{displaymath}
        \xymatrix{{\bf e} : & 0\ar[r] & \Ker f\ar[r] &
M(\beta\sigma)\ar[r]^{f} & M(\eta\alpha)\ar[r] & \Coker f\ar[r] & 0
        }        
\end{displaymath}
represents a non-zero element of $\Ext{2}{C}(\Coker f, \Ker f)$ if and
only if $\alpha\beta$ lies in $I$.
\end{lemma}
{\bf Proof.}
In view of \ref{lemm::HRS}, the sequence ${\bf e}$
represents a non-zero element of $\Ext{2}{C}(\Coker f, \Ker f)$ if and
only if there exists no short exact sequence of the form
\begin{displaymath}
        \xymatrix{0\ar[r] & M(\beta\sigma)\ar[r]^{(p,
g)^{t}\phantom{xx}} & \Ima f \oplus X \ar[r]^{\phantom{xx}(j,h)} &
M(\eta\alpha)\ar[r] & 0,}
\end{displaymath}
where $f=jp$ is the canonical factorisation.  

Assume such a sequence exists.  Since $S_b$ appears exactly once as a
composition factor of $M(\beta\sigma)$ and $M(\eta\alpha)$, then it also
appears exactly once as a composition factor of $X$.  Therefore, there
exists a unique indecomposable summand $Y$ of $X$ admitting $S_b$ as a
composition factor.  

We claim that  $g:M(\beta\sigma)\longrightarrow X$ is a monomorphism: let indeed $x \in Q_0$ and take a vector $v \in M(\beta\sigma)_x$ such that $g_x(v)=0$. 
If $x \neq b$, then $p_x(v)=0$ and so $(p,g)^t_x(v)=0$ which implies $v=0$.
If $x=b$, then $(p_a(\beta v), g_a(M(\beta\sigma)_{\beta}(v)))^t = (p,g)^t_a(M(\beta\sigma)_{\beta}(v))=(\Ima f \oplus X)_{\beta}(p,g)^t_b(v)=0$ which implies $M(\beta\sigma)_{\beta}(v)=0$. Since $(M(\beta\sigma)_{\beta}$ is injective, then $v=0$. This completes the proof of our claim.

Since the evaluation $M(\beta\sigma)_{\beta}$
of the module $M(\beta\sigma)$ on the arrow $\beta$ is non-zero, we must
have $X_{\beta} \neq 0$.  Now, $S_b$ is a composition factor of $Y$,
hence $Y_{\beta} \neq 0$ as well.  Similarly, $h$ is an epimorphism and it follows that 
$Y_{\alpha} \neq 0$.
On the other hand, $Y$ must be a string or a band module.  The above
reasoning implies that $\alpha\beta$ must then be a subpath of a string
or a band, which implies that $\alpha\beta \notin I$, as required.

Conversely, if $\alpha\beta \notin I$, then we have a short exact sequence
\begin{displaymath}
        \xymatrix{ 0\ar[r] & M(\beta\sigma)\ar[r] & S_{b} \oplus
M(\eta\alpha\beta\sigma) \ar[r] & M(\eta\alpha)\ar[r] & 0,
        }
\end{displaymath}
and hence ${\bf e}$ represents the zero element in
$\Ext{2}{C}(\Coker f, \Ker f)$.
\hfill $\Box$

\subsubsection{Arrows}

From now on, let $C$ be a tilted algebra of type $\mathbb{A}$ or
$\tilde{\mathbb{A}}$.  
We  give a description of the elements of $\tilde{C} = C \ltimes
\Ext{2}{C}(DC, C)$ corresponding to the arrows of its ordinary quiver.
 In \cite[(2.4)]{ABS}, it is proved that the quiver of $\tilde{C}$ is
obtained from that of $C$ by adding an arrow from $x$ to $y$ for each
relation from $y$ to $x$.  The elements of $\tilde{C}$ corresponding to
the arrows of $C$ are of the form $(\alpha, 0)$, where $\alpha$ is an
arrow of $C$.  

The other arrows correspond to relations in $C$.  Let $\alpha\beta$  be
a relation from $c$ to $a$ in $C$, and let $\xi_{\alpha\beta}$ be the
corresponding new arrow in $\tilde{C}$.  

\begin{lemma}
The new arrow $\xi_{\alpha\beta}$ lies in $0 \oplus \Ext{2}{C}(I_c, P_a)$.
\end{lemma}
{\bf Proof.} This new arrow lies in $e_a \tilde{C} e_c$, which can be written
as the direct sum $e_a C e_c \oplus e_a \textrm{Ext}^2_C(DC,C) e_c$.  We
know from (\ref{prop::bound quiver}) that the quiver of $C$ contains no
double-zero.  Consequently, there are no paths from $a$ to $c$, and
hence $e_a C e_c = 0$.  Moreover, $e_a \textrm{Ext}^2_C(DC,C) e_c =
\textrm{Ext}^2_C(I_c, P_a)$.  The element $\xi_{\alpha\beta}$ thus lies
in $0 \oplus \textrm{Ext}^2_C(I_c, P_a)$.
\hfill $\Box$  
\bigskip

The following lemma gives the dimension and a basis of the extension
space involved in the last expression.

\begin{lemma}\label{lemm::fleches}
Let $\alpha : c \longrightarrow b$ and $\beta : b \longrightarrow a$ be
two arrows of $C$ such that $\alpha\beta \in I_C$.  
\begin{description}
\item[(a)] The dimension of the vector space $\Ext{2}{C}(I_c, P_a)$ is 1
or 2.  Moreover, its dimension is 2 if, and only if, the following
situation occurs in the bound quiver of  $C$ :
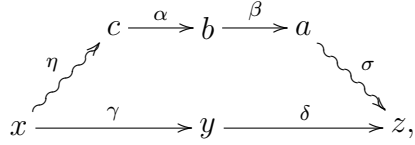
\begin{displaymath}
        \xymatrix{ & c\ar[r]^{\alpha} & b\ar[r]^{\beta}  &
a\ar@{~>}[dr]^{\sigma} & \\
           x\ar[rr]^{\gamma}\ar@{~>}[ur]^{\eta} &        &
y\ar[rr]^{\delta} &   &z        ,
        }
\end{displaymath}
where $\gamma$ and $\delta$ are arrows, $\eta$ and $\sigma$ are paths,
possibly stationary, without relations, and $\alpha\beta, \gamma\delta$
are relations.

\item[(b)] If the dimension of the space is 1, then a basis is given by
the sequence
\begin{displaymath}
 \xymatrix{ {\bf e}_1: & 0 \ar[r] & P_a \ar[r] &
M(\beta\sigma)\ar[r] & M(\eta\alpha)\ar[r] & I_c\ar[r] & 0,
 }        
\end{displaymath}
where $\eta$ and $\sigma$ are paths such that $I_c = M(\eta)$ and $P_a =
M(\sigma)$.

\item[(c)] If the dimension of the space is 2, then a basis is given by
the sequences
\begin{displaymath}
 \xymatrix{ {\bf e}_1: & 0 \ar[r] & P_a \ar[r] &
M(\beta\sigma)\ar[r] & M(\eta\alpha)\ar[r] & I_c\ar[r] & 0
 }        
\end{displaymath}
\begin{displaymath}
 \xymatrix{ {\bf e}_2: & 0 \ar[r] & P_a \ar[r] &
M(\sigma\delta^{-1})\ar[r] & M(\gamma^{-1}\eta)\ar[r] & I_c\ar[r] & 0,
 }        
\end{displaymath}
where $\gamma$, $\delta$, $\eta$ and $\sigma$ are as in the figure in (a).
\end{description}
\end{lemma}
\bigskip\bigskip

{\bf Proof.}
\begin{description}
\item[(a)] It is known from \cite{ABS} that there is a new arrow from
$a$ to $c$; thus the dimension cannot be zero.  On the other hand, since
$C$ is gentle and without double-zero, the local situation of the relation
$\alpha\beta$  can be described by the following figure, where dotted
lines represent relations.
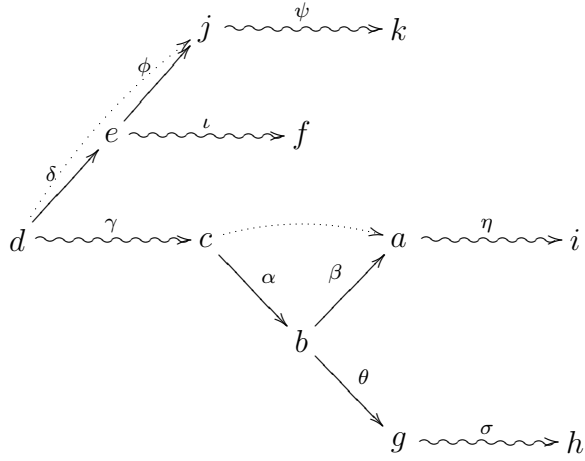
\begin{displaymath}
        \xymatrix{ &  & j\ar@{~>}[rr]^{\psi} & & k &  &\\
                   & e\ar@{~>}[rr]^{\iota}\ar[ur]^{\phi} & & f & &        \\              
d\ar@{~>}[rr]^{\gamma}\ar[ur]^{\delta}\ar@/^/@{.>}[uurr] & &
c\ar[dr]^{\alpha}\ar@/^/@{.>}[rr] & & a\ar@{~>}[rr]^{\eta} & & i\\
                    & & & b\ar[ur]^{\beta}\ar[dr]^{\theta} & & & \\
                    & & & & g\ar@{~>}[rr]^{\sigma} & & h
        }
\end{displaymath}
This diagram allows us to compute a projective resolution of $I_c$ in \\
mod $C$ :
\begin{displaymath}
        \xymatrix{ 0 \ar[r]^{p_3} & P(2)\ar[r]^{p_2} & P(1)\ar[r]^{p_1}
&  P(0)\ar[r]^{p_0} & I_c\ar[r] & 0.
        }
\end{displaymath}
where $P(2) = M(\psi) \oplus M(\eta)$, $ P(1) = M(\iota^{-1}\phi\psi)
\oplus M(\sigma^{-1}\theta^{-1}\beta\eta)$ and
$P(0)=M(\iota^{-1}\delta^{-1}\gamma\alpha\theta\sigma)$.
Note that some direct summands of the terms of this sequence can be
zero.  Applying $\Hom{C}(-,P_a)$, we get a complex
\begin{displaymath}
        \xymatrix{ 0 \ar[r] & \Hom{C}(I_c,P_a)\ar[r]^{(p_0,P_a)} &
\Hom{C}(P(0), P_a)\ar[r]^{\phantom{xxxxxxxx}(p_1,P_a)} &   \\
        \ar[r]^{(p_1,P_a)\phantom{xxxxxxxx}} & \Hom{C}(P(1),
P_a)\ar[r]^{(p_2,P_a)} & \Hom{C}(P(2),
P_a)\ar[r]^{\phantom{xxxxxxxx}(p_3,P_a)} & 0.
        }.
\end{displaymath}
This yields
\begin{displaymath}
        \Ext{2}{C}(I_c, P_a) = \frac{\Ker \Hom{}(p_3, P_a)}{\Ima
\Hom{}(p_2, P_a)} = \frac{\Hom{}(M(\psi), P_a) \oplus \Hom{}(M(\eta),
P_a)}{\Ima \Hom{}(p_2, P_a)}.
\end{displaymath}

Since $P_a = M(\eta)$, then $\dim \Hom{}(M(\eta), P_a) = 1$, and since
\\ $\Hom{}(M(\sigma^{-1}\theta^{-1}\beta\eta), P_a) = 0$, no non-zero morphism in
$\Hom{}(M(\eta), P_a)$ factors through $p_2$.

Moreover, we claim that $\Hom{}(M(\psi), P_a)$ is non-zero if and only if $j = i$.
Indeed, a non-zero morphism from $M(\psi)$ to $P_a$ can only exist when $j$ coincides with a vertex on the path $\eta$. 
But if $j$ were a vertex different from $i$, then there would be an arrow $\phi': j \to j'$ in the path $\eta$, forcing the relation $\phi\phi''$ and creating the double-zero $\delta\phi\phi'$.
Thus $j = i$.
 In this case, $\psi$ has no choice but to be the trivial path in $i$,
and $\dim \Hom{}(M(\psi), P_a) = 1$.  Since
$\Hom{}(M(\iota^{-1}\phi\psi, P_a) = 0$, no non-zero morphism in $\Hom{}(M(\psi),
P_a)$ factors through $p_2$.  

Hence no non-zero morphism in $\Hom{}(M(\psi), P_a) \oplus \Hom{}(M(\eta), P_a)$
factors through $p_2$.  Thus the dimension of this space is either $1$
or $2$, and it is $2$ exactly when $i=j$.  In this case, and in this
case only, we have  
\begin{displaymath}
        \xymatrix{ & c\ar[r]^{\alpha} & b\ar[r]^{\beta}  & a\ar@{~>}[dr]& \\
              d\ar[rr]\ar@{~>}[ur] &        & e\ar[rr] &   & j        
                     }
\end{displaymath}
as desired.

\item[(b)] It follows from (\ref{lemm::zero}) that ${\bf e}_1$ is
non-zero.  The result follows.

\item[(c)] It follows from (\ref{lemm::zero}) that ${\bf e}_1$ and
${\bf e}_2$ are non-zero.  

It remains to be shown that ${\bf e}_1$ and ${\bf e}_2$ are
linearly independant.  Suppose there exists a non-zero scalar $\lambda$
such that ${\bf e}_2 + \lambda{\bf e}_1 = 0$.
Computing this sum, we get the sequence
\begin{displaymath}
        \xymatrix{ 0\ar[r] & P_a\ar[r] & M(\beta\sigma\delta^{-1})\ar[r]^f & M(\gamma^{-1}\eta\alpha)\ar[r] & I_c\ar[r] &0
        }
\end{displaymath}
where all morphisms are multiples of the natural morphisms between
string modules.

Here, applying (\ref{lemm::zero}) is not possible, since $\Ima f = S_b
\oplus S_y$, but a similar technique of proof can be used.

Suppose there exist a module $X$ and morphisms $g$ and $h$ such that
the sequence 
\begin{displaymath}
        \xymatrix{0\ar[r] & M(\beta\sigma\delta^{-1})\ar[r]^{(p,g)^{t}\phantom{xx}} & (S_b \oplus S_y) \oplus X\ar[r]^{\phantom{xx}(j,h)} & M(\gamma^{-1}\eta\alpha) \ar[r] & 0
        }           
\end{displaymath}  
is exact, where $f = jp$ is the canonical factorisation.  Since $S_b$
appears exactly once as a composition factor of
$M(\beta\sigma\delta^{-1})$ and $M(\gamma^{-1}\eta\alpha)$, then it also
appears exactly once as a composition factor of $X$.  Therefore, there
exists a unique indecomposable summand $Y$ of $X$ admitting $S_b$ as a
composition factor. As in the proof of \ref{lemm::zero}, we show that
$Y_{\beta} \neq 0$ and $Y_{\alpha} \neq 0$.  Therefore, $\alpha\beta$
must be a subpath of a string or a band, which is a contradiction, since
it is a relation.

The sequences ${\bf e}_1$ and ${\bf e}_2$ thus form a basis
of the extension space.\hfill $\Box$  
\end{description}

It remains to determine which of the basis elements are represented by  arrows of
$\tilde{C}$.

\begin{lemma}
Let $\alpha : c \longrightarrow b$ and $\beta : b \longrightarrow a$ be
two arrows of the quiver of $C$ such that $\alpha\beta$ is a relation.
 Let $\xi_{\alpha\beta}$ be the corresponding new arrow in $\tilde{C}$.
 With the notation of (\ref{lemm::fleches}), the element $(0,
{\bf e}_1)$ can be chosen to represent $\xi_{\alpha\beta}$.
\end{lemma}
{\bf Proof.} The space $0 \oplus \Ext{2}{C}(I_c, P_a)$ contains at least one
arrow.  

If its dimension is 1, then the result is obvious.

If its dimension is 2, then lemma (\ref{lemm::fleches}) describes the
situation of $\alpha\beta$ in the quiver of $C$.  Two cases arise.

First,  suppose that $\eta$ and $\sigma$ are both trivial paths.
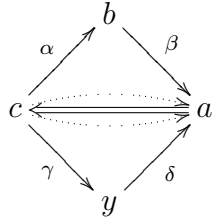
\begin{displaymath}
        \xymatrix{ & b\ar[dr]^{\beta} & \\
                   c\ar[ur]^{\alpha}\ar@/^/@{.>}[rr] \ar@/_/@{.>}[rr]
\ar[dr]_{\gamma} & & a\ar@{=>}[ll] \\
                    & y\ar[ur]_{\delta} &
        }
\end{displaymath}
In this case, two arrows from $a$ to $c$ are added to the quiver.  Both
$(0, {\bf e}_1)$ and $(0, {\bf e}_2)$ must thus represent
arrows of $\tilde{C}$.

Second, suppose $\eta$ and $\sigma$ are not both trivial.  In this case,
lemma (\ref{lemm::fleches}) implies that  $\Ext{2}{C}(I_x, P_z)$ is of
dimension 1, and that a basis is given by
\begin{displaymath}
        \xymatrix{ {\bf e'} : & 0\ar[r] & P_z\ar[r] &
M(\delta)\ar[r] & M(\gamma)\ar[r] & I_x\ar[r] & 0.
        }
\end{displaymath}
Reasoning as above, we get that $(0, {\bf e'})$ represents the new
arrow from $z$ to $x$.  Moreover, a straightforward calculation yields
$(\sigma, 0)(0, {\bf e'})(\eta, 0) = (0, {\bf e}_2)$.

Since one of $\eta$ and $\sigma$ is not trivial, one of $(0, \eta)$ and
$(0, \sigma)$ must lie in $\textrm{rad } \tilde{C}$.  Therefore $(0,
{\bf e}_2) \in \textrm{rad}^2 \tilde{C}$, and $(0,
{\bf e}_1) \in \textrm{rad } \tilde{C} \setminus \textrm{rad}^2
\tilde{C}$; in other words, $(0, {\bf e}_1)$ represents an arrow
from $a$ to $c$.
\hfill $\Box$  
\bigskip

\subsubsection{Relations}
Knowing how to write arrows in $\tilde{C}$ allows us to compute the
relations.

\begin{lemma}\label{lemm::relations}
Let $C=kQ_C/I_C$ and $\tilde{C} = kQ_{\tilde{C}}/I_{\tilde{C}}$.

\begin{enumerate}

\item Let $\omega_1, \omega_2, \ldots, \omega_n$ be paths from $x$ to
$y$ in the quiver of $C$, and let $\lambda_1, \lambda_2, \ldots,
\lambda_n \in k$.  Then $\sum_{i=1}^n \lambda_i (\omega_i, 0) = 0$ in
$\tilde{C}$ if, and only if, $\sum_{i=1}^n \lambda_i \omega_i = 0$ in $C$.

\item Let $\alpha : c \longrightarrow b$ and $\beta : b \longrightarrow
a$ be two arrows in the quiver of $C$ such that $\alpha\beta$ is a
relation.  Let $(0, {\bf e}_1)$ be the element representing the
corresponding new arrow, where ${\bf e}_1$ is as in lemma
(\ref{lemm::fleches}).  Then $(0, {\bf e}_1)(\alpha, 0) = 0$ and
$(\beta, 0)(0, {\bf e}_1) = 0$.

\item The ideal $I_{\tilde{C}}$ is generated by the relations of $C$ and
those described in 2.

\end{enumerate}
\end{lemma}
{\bf Proof.} \begin{enumerate}

\item This is shown in \cite{ACT}.

\item Viewing $\alpha$ as an element of $\textrm{End }DC$, more
precisely as a morphism from $I_b$ to $I_c$, we can compute
${\bf e}_1\beta$:

\begin{displaymath}
        \xymatrix{ {\bf e}_1\beta : & 0\ar[r] & P_a\ar[r] &
M(\beta\sigma)\ar[r] & M(\eta\alpha) \oplus M(\varphi\gamma)\ar[r] &
I_b\ar[r] & 0,
        }
\end{displaymath}
where $I_b = M(\eta\alpha\gamma^{-1}\varphi^{-1})$.  This sequence
represents the zero element, because of (\ref{lemm::HRS}) and exactness of the sequence
{\small
\begin{displaymath}
        \xymatrix{ 0 \ar[r] & M(\beta\sigma)\ar[r] & S_b \oplus
M(\varphi\gamma\beta\sigma) \oplus M(\eta\alpha)  \ar[r] &
M(\varphi\gamma)\oplus M(\eta\alpha) \ar[r] & 0.
        }
\end{displaymath}
}

Therefore $(0, {\bf e}_1)(\alpha, 0) = 0$.  

In a dual way, we prove that $(\beta, 0)(0, {\bf e}_1) = 0$.

\item It is sufficient to show that new arrows in the quiver of
$\tilde{C}$ are not involved in other relations than those described in 2.  

First suppose that $w$ is a monomial relation involving new arrows and
other than those relations described in 2.  Then it must contain exactly one new arrow
$\xi$, corresponding to a relation $\alpha\beta$; otherwise the quiver
of $C$ would contain a double-zero.  Write $w = u\xi v$, where $u$ and
$v$ are non-zero paths consisting of arrows of $C$.  Let
${\bf e}_1$ be the sequence as in (\ref{lemm::fleches})
corresponding to $\xi$.  Then $(u,0)(0,{\bf e}_1)(v,0) = (0,
u{\bf e}_1 v)$, where $u{\bf e}_1 v$ is the sequence
{\small
\begin{displaymath}
        \xymatrix{ 0\ar[r] & M(u^{-1}u')\ar[r] & M(\beta u^{-1}u')\ar[r]
& M(v'v^{-1}\alpha)\ar[r] & M(v'v^{-1})\ar[r] & 0,
        }
\end{displaymath}
}

where $u'$ and $v'$ are paths in the quiver of $C$.  The figure below
illustrates the local situation, where $\alpha\beta = \gamma'\delta' =
0$, the last arrow of $u$ and the first of $t$ form a relation, as for
the last of $v$ and $v'$ and the first of $w'$ and $w$, respectively.
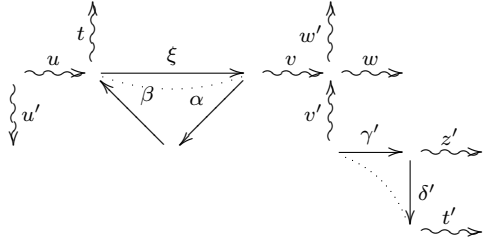
\begin{displaymath}
        \xymatrix{ & & & & & & \\
                   \ar@{~>}[r]^{u}\ar@{~>}[d]^{u'} &
\ar@{~>}[u]^{t}\ar[rr]^{\xi} & & \ar[dl]_{\alpha}\ar@/^/@{.}[ll]
 \ar@{~>}[r]^{v} & \ar@{~>}[r]^{w}\ar@{~>}[u]^{w'} & & \\
                   & & \ar[ul]_{\beta} & &
\ar@{~>}[u]^{v'}\ar[r]^{\gamma'}\ar@/^/@{.}[dr] &
\ar@{~>}[r]^{z'}\ar[d]^{\delta'} & \\
                   & & & & & \ar@{~>}[r]^{t'} &
        }
\end{displaymath}

This yields the following commutative diagram, where the first line is a
projective resolution of $M(v'v^{-1})$ :
{\small
\begin{displaymath}
        \xymatrix{ 0\ar[r] & P(2) \ar[r]^{p_2}\ar[d]^{(f,0)} & P(1)
\ar[r]^{p_1}\ar[d]^{(0,g,0)} & P(0) \ar[r]^{p_0}\ar[d]^{(h,\ell)} &
M(v'v^{-1}) \ar[r]\ar@{=}[d] & 0 \\
                   0\ar[r] & M(u^{-1}u')\ar[r] & M(\beta u^{-1}u')\ar[r]
& M(v'v^{-1}\alpha)\ar[r] & M(v'v^{-1})\ar[r] & 0,
        }
\end{displaymath}
}

where $P(2) = M(t)\oplus M(t')$, $P(1) = M(w'w)\oplus
M(z^{-1}\gamma^{-1}\beta t)\oplus M(z'^{-1}\delta't')$ and $P(0) =
M(z^{-1}\gamma^{-1}\alpha^{-1}vw)\oplus M(z'^{-1}\gamma'^{-1}v'w')$, and
all non-zero morphisms are the natural morphisms between string modules.
 It is then seen that $(f,0)$ cannot factor through $p_2$, and thus the
lower exact sequence is non-zero.
  Hence there are no other monomial relations than those in 2.

Now suppose we have a minimal relation of the form $\sum_{i=1}^{m}
\lambda_i w_i$, where each $\lambda_i$ is a non-zero scalar, each $w_i$
is a path in the quiver of $\tilde{C}$ and $m \geq 2$.  At least one of
the $w_i$ must pass through a new arrow, and since $C$ contains no
double zero, this implies that each $w_i$ must pass through exactly one
new arrow, say $\xi_i$, corresponding to a relation $\alpha_i \beta_i$.
 Write $w_i = u_i \xi_i v_i$, where $u_i$ and $v_i$ are paths of the
quiver of $C$.

Since the quiver of $C$ contains at most one cycle, we must have $m =
2$.  Since $k$ is a field, we may suppose that $\lambda_1  = 1$.
 Letting ${\bf e}_1$ and ${\bf e}_2$ be the sequences
associated to $\xi_1$ and $\xi_2$, respectively, we get that $u_1
{\bf e}_1 v_1$ and $\lambda_2 u_2 {\bf e}_2 v_2$ are both
sequences of the form above.  Their sum is the sequence

\begin{displaymath}
        \xymatrix{ 0 \ar[r] & M(u_2^{-1}u_1)\ar[r] &
M(\beta_2u_2^{-1}u_1\beta_1^{-1}) \ar[r] &  \\
         \ar[r] & M(\alpha_2^{-1}v_2v_1^{-1}\alpha_1) \ar[r] &
M(v_1v_2^{-1}) \ar[r] & 0.
        }
\end{displaymath}

By an argument similar to the one given in the proof of
(\ref{lemm::fleches})(c), this element is not zero, a contradiction.
 Hence no binomial relations exist in $\tilde{C}$.
\end{enumerate}
\hfill $\Box$  
\bigskip

The relations described in the preceding lemma make $\tilde{C}$ a gentle
algebra.  
\begin{lemma}\label{lemm::derniere implication}
If $\tilde{C}$  is cluster-tilted of type $\mathbb{A}$ or
$\tilde{\mathbb{A}}$, then $\tilde{C}$ is gentle.
\end{lemma}
{\bf Proof.} The relations of $\tilde{C}$ are known (see
\ref{lemm::relations}).  Moreover, $C$ is gentle.  

Suppose that there are $r$ new arrows.  Let us add the new arrows and
the corresponding new relations one by one, thus obtaining a sequence $C
= C_0 , C_1, \ldots, C_r = \tilde{C}$ of algebras.  We show that $C_i$
is gentle for all $i$ in $\{ 0,1,2,\ldots,r \}$.

Since $C$ is gentle, then so is $C_0$.  Suppose that $C_i$ is gentle,
where $i$ is in $\{ 0,1,2,\ldots,r-1 \}$.  To get $C_{i+1}$, we add one
new arrow, say $\gamma$ from $x$ to $y$.  This arrow comes from a
relation $\alpha\beta$ from $y$ to $x$ in $C$.  We must add the
relations $\beta\gamma$ and $\gamma\alpha$ to obtain $C_{i+1}$.

Since $C_i$ is gentle, there were already at most two  arrows starting
from $x$ in $C_i$.  Suppose that there were two, say $\eta_1$ and
$\eta_2$.  Since $C_i$ is gentle, then $\beta$ is involved in a relation
with one of the two, say $\eta_1$.  The arrow $\eta_1$ cannot be in $C$,
otherwise there would be a double zero involving $\alpha\beta$ and
$\beta\eta_1$.  So the arrow $\eta_1$ comes from a relation
$\sigma\beta$ in $C$.  Since $C$ is gentle, we must have that $\sigma =
\alpha$, so that $\eta_1 = \gamma$, which is absurd because $\gamma$ is
not in $C_i$.  

Therefore, in $C_i$, there is at most one outgoing arrow from $x$, and
this arrow is not involved in a relation with $\beta$.  This shows that
in $C_{i+1}$, there are at most two arrows starting from $x$, say $\eta$
and $\gamma$, and that $\beta\eta$ is not a relation while $\beta\gamma$
is.  Moreover, there is at most one more  arrow ending in $x$, say
$\delta$, and since $C_i$ is gentle, we have that $\delta\eta$ is a
relation, while $\delta\gamma$ is not.  So the relations at $x$ are
those found in a gentle algebra.

Using a similar argument for the vertex $y$, we get that $C_{i+1}$ is a
gentle algebra.

By induction, $\tilde{C}$ is a gentle algebra.
\hfill $\Box$  
\bigskip

\begin{example} \rm
Lemma \ref{lemm::relations} allows us to compute the relation-extension
of any gentle tilted algebra.  As an illustration, consider the two
algebras given in example \ref{exem::an et antilde}.  The
relation-extension of each is given in the following diagram:
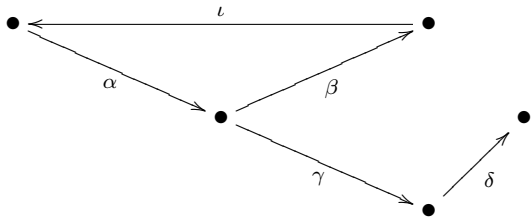
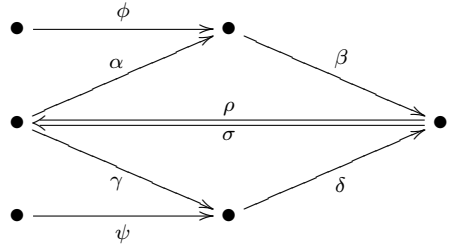
\begin{displaymath}
        \xymatrix{ \bullet\ar[dr]_{\alpha} & & \bullet\ar[ll]_{\iota} &
                  & \bullet\ar[r]^{\phi} & \bullet\ar[dr]^{\beta} &   \\
                   & \bullet\ar[ur]_{\beta}\ar[dr]_{\gamma} & & \bullet
   & \bullet\ar[ur]^{\alpha}\ar[dr]_{\gamma} & &
\bullet\ar@{=>}[ll]^{\sigma}_{\rho} \\
                   & & \bullet\ar[ur]_{\delta} &                        
  & \bullet\ar[r]_{\psi} & \bullet\ar[ur]_{\delta} &  \\
                   & \alpha\beta = \iota\alpha = \beta\iota = 0 & &    
                              & & \alpha\beta = \gamma\delta =
\rho\alpha = 0 \\
                   & & &                                                
  & & \beta\rho = \sigma\gamma = \delta\sigma =0 &.
        }
\end{displaymath}
\end{example}

%
%
\subsubsection{Proof of the main theorem}\label{sect::preuve principale}
 Now the proof of theorem \ref{theo::principal}, which has been written
in separate parts during the last sections,  can be stated properly.
\bigskip

{\bf Proof (of theorem \ref{theo::principal})}:
1 implies 2 is shown in  (\ref{lemm::1}),   2 implies 4 in  (\ref{lemm::3}), 
4 implies 3 in (\ref{lemm::derniere implication}) and 3 implies 1 in (\ref{lemm::2}).
\hfill $\Box$  
\bigskip

%
%

\section{Geometry of surfaces and $A(\Gamma)$}
We study in this section more connections between geometric 
properties of the marked surface $(S,M)$ and properties of the 
algebra $A(\Gamma)$ given by a triangulation of $(S,M)$.

\subsection{Cluster-tilted algebras arising from surfaces}

We first address the question which of the algebras $A(\Gamma)$  
are cluster-tilted. Recall that all algebras $A(\Gamma)$  share the 
properties (2) and (3) from Theorem \ref{structure} with every 
 cluster-tilted algebra.
Moreover, it is shown in \cite{CCS} and   \cite{BV} that the cluster-tilted 
algebras of type $\mA$ are algebras $A(\Gamma)$ arising from 
a triangulation of an unpunctured polygon.
In this section, we show the following generalization:

\begin{theorem}
Let $A(\Gamma)$ be the algebra associated to  the triangulation
$\Gamma$ of an unpunctured marked surface $(S,M)$.
 Then the following statements are equivalent:                                                                                                                     
\begin{itemize}
\item[(1)] the algebra $A(\Gamma)$ is cluster-tilted

\item[(2)] the algebra $A(\Gamma)$ is cluster-tilted of type $\mathbb{A}$ or
$\tilde{\mathbb{A}}$
\item[(3)]  $S$ is a disc or an annulus

\end{itemize}             
Moreover, all cluster-tilted algebras of type  $\mathbb{A}$ (or 
$\tilde{\mathbb{A}}$) are of the form $A(\Gamma)$ for some triangulation $\Gamma$
of a disc $S$ (or an annulus $S$, respectively).
                                                                                         
\end{theorem}

{\bf Proof.}  
It is clear that $(2)$ implies $(1)$. Let us show the
converse: Suppose that the algebra $A(\Gamma)$ is cluster-tilted. Thus
there is a sequence of mutations transforming the quiver with potential
defining $A(\Gamma)$ into some quiver $Q$ with zero potential. This
sequence of mutations corresponds to a sequence of flips, transforming
the triangulation $\Gamma$ of $(S,M)$ into a triangulation $T$ with
$Q(T)=Q$ and zero potential. Hence $A(T)=kQ$ is hereditary. Since we
 know from  (\ref{structure})
 that $A(T)$ is gentle, this leaves only the possibilities
that $Q$ is of type $\mA$ or
$\widetilde{\mA}$. Therefore the algebra $A(\Gamma)$ is cluster-tilted 
of type $\mathbb{A}$ or $\tilde{\mathbb{A}}$.
\bigskip

We prove now the equivalence of $(2)$ and $(3)$. Since all
triangulations on $(S,M)$ are flip-equivalent (see \cite{Ha}) and flips of the
triangulation correspond to mutations of the corresponding quiver with
potential (see \cite{La}), it is sufficient to consider one particular triangulation.
In case $S$ is a disc, we choose the triangulation to  be in the form of
a fan, giving rise to a linear oriented quiver of type $\mA$.
In case where $S$ is an annulus, we choose the triangulation given by
two fans in opposite direction as shown in the following figure
(identify the left and right vertical edge):

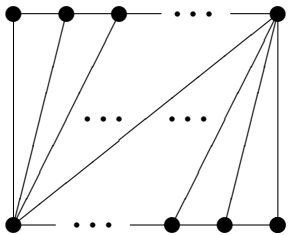
\begin{center}
\setlength{\unitlength}{2pt}
\begin{picture}(50,40)(0,0)
   \put(0,0){\vertex}
   \put(30,0){\vertex}
   \put(40,0){\vertex}
   \put(50,0){\vertex}
   \put(0,40){\vertex}
   \put(10,40){\vertex}
   \put(20,40){\vertex}
   \put(31,40){\svertex}
   \put(34,40){\svertex}
   \put(37,40){\svertex}
   \put(14,20){\svertex}
   \put(17,20){\svertex}
   \put(20,20){\svertex}
  \put(12,0){\svertex}
   \put(15,0){\svertex}
   \put(18,0){\svertex}
   \put(30,20){\svertex}
   \put(33,20){\svertex}
   \put(36,20){\svertex}
   \put(50,40){\vertex}
   \put(50,0){\vertex}
   \put(30,0){\line(1,2){20}}
   \put(40,0){\line(1,4){10}}
 \put(0,0){\line(5,4){50}}
   \put(0,0){\line(1,2){20}}
  \put(0,0){\line(1,4){10}}
   \put(0,40){\line(1,0){28}}
   \put(50,0){\line(-1,0){28}}
   \put(41,40){\line(1,0){8}}
   \put(0,0){\line(1,0){8}}
   \put(0,0){\line(0,1){40}}
   \put(50,0){\line(0,1){40}}
  \end{picture}
\end{center}
\bigskip

The corresponding quiver is of type $\widetilde{\mA}$ with zero potential, 
thus (3) implies (2).
Conversely, we know from Proposition \ref{topology} that the quivers $Q(\Gamma)$ uniquely determine the topology of the unpunctured marked surface $(S,M)$. Therefore $S$ is a disc or an annulus, respectively, and since all triangulations are flip-equivalent, it is clear that all cluster-tilted algebras of the corresponding type occur.
\hfill $\Box$

\subsection{Curves in $(S,M)$ and string modules}

In this section we are comparing strings in $A(\Gamma)$ to curves in $(S,M)$. By a curve in $(S,M)$ we mean a curve $\gamma$ in $S$ whose endpoints lie in $M$ and where all points except the endpoints lie in the interior of $S$. We usually consider curves up to homotopy. For instance,  for two distinct curves $\gamma$ and $\delta$ in $(S,M)$, the intersection number $I_{\Gamma}(\gamma,\delta)$ is defined as the minimal number of transversal intersections of two representatives of the homotopy classes of $\gamma$ and $\delta$.
Denote the internal arcs of the triangulation $\Gamma$ by $\{ a_1, \ldots, a_n \}$.
Then we define the {\em intersection vector} $ I_{\Gamma}(\gamma)$ of a curve $\gamma$ as 
$$ I_{\Gamma}(\gamma) = ( I_{\Gamma}(\gamma, a_1) \, \ldots , I_{\Gamma}(\gamma,a_n))$$

\begin{proposition}
Let $\Gamma$ be a triangulation of an unpunctured marked  surface $(S,M)$.  Then
there exists a bijection $\{\gamma\} \mapsto w(\gamma)$ between the homotopy classes of curves in $(S,M)$ not homotopic to an arc in $\Gamma$ and the strings of $A(\Gamma)$.  Under this bijection, the intersection vector corresponds to the dimension vector of the corresponding string module, that is 
$$I_{\Gamma}(\gamma) = \udim M(w(\gamma))$$
\end{proposition}

{\bf Proof.} 
Let $ w = \;x_1 \stackrel{\alpha_1}{\longleftrightarrow} x_2 \stackrel{\alpha_2}{\longleftrightarrow}\cdots  \stackrel{\alpha_{s-1}}{\longleftrightarrow} x_s$ be a string in $A(\Gamma)$.
We define a curve $\gamma(w)$ in $(S,M)$ as follows: The arcs $x_1$ and $x_2$ belong to the same triangle $T_1$ since they are joined by an arrow in $A(\Gamma)$. We connect the midpoints of $x_1$ and $x_2$ by a curve $\gamma_1$ in the interior of $T_1$. 
Proceeding in the same way with the remaining arcs $x_2, \ldots, x_s$ we obtain curves $\gamma_2, \ldots, \gamma_{s-1}$ connecting the midpoints of the respective arcs. 
The internal arc $x_1$ belongs to two triangles, the triangle $T_1$ which we considered above and another triangle $T_0$.
Let $P \in M$ be the marked point in $T_0$ opposite to the arc $x_1$.
We now connect $P$ with the midpoint of $x_1$ by a curve $\gamma_0$ in the interior of $T_0$, and proceed in the same way on the other end of the string $w$, connecting the midpoint of $x_s$ with a marked point $Q$ by some curve $\gamma_s$.
The curve $\gamma(w)$ is then defined  as the concatenation of the curves $\gamma_0, \ldots, \gamma_s$.
\bigskip

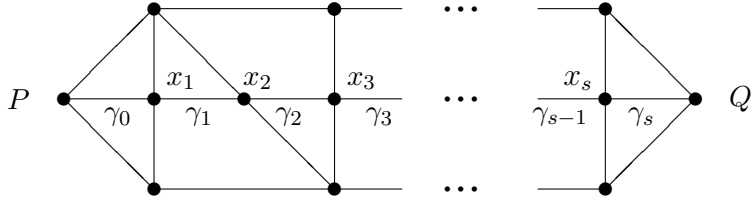
\setlength{\unitlength}{0.6mm}
\begin{center}
\begin{picture}(140,40)(0,0)
   \put(0,20){\vertex}
   \put(20,0){\vertex}
   \put(20,40){\vertex}
   \put(40,20){\vertex}
   \put(60,20){\vertex}
   \put(120,20){\vertex}
   \put(120,0){\vertex}
   \put(120,40){\vertex}
   \put(140,20){\vertex}
   \put(60,40){\vertex}
   \put(20,20){\vertex}
   \put(60,0){\vertex}
   \put(60,40){\vertex}

  \put(91,0){\svertex}
   \put(85,0){\svertex}
   \put(88,0){\svertex}
  \put(91,20){\svertex}
   \put(85,20){\svertex}
   \put(88,20){\svertex}
  \put(91,40){\svertex}
   \put(85,40){\svertex}
   \put(88,40){\svertex}
   
   \put(0,20){\line(1,1){20}}
   \put(0,20){\line(1,-1){20}}
   \put(0,20){\line(1,0){75}}
   \put(140,20){\line(-1,0){35}}
   \put(140,20){\line(-1,1){20}}
   \put(140,20){\line(-1,-1){20}}
   \put(20,40){\line(1,-1){40}}
   \put(20,40){\line(1,0){55}}
   \put(20,0){\line(1,0){55}}
   \put(20,40){\line(0,-1){40}}
   \put(120,40){\line(0,-1){40}}
   \put(120,40){\line(-1,0){15}}
   \put(120,0){\line(-1,0){15}}
   \put(60,40){\line(0,-1){40}}

  \put(-10,20){\makebox(0,0){\small $P$}}
  \put(150,20){\makebox(0,0){\small $Q$}}
   \put(26,24){\makebox(0,0){\small $x_1$}}
   \put(43,24){\makebox(0,0){\small $x_2$}}
   \put(66,24){\makebox(0,0){\small $x_3$}}
   \put(114,24){\makebox(0,0){$x_s$}}
    \put(30,16){\makebox(0,0){\small $\gamma_1$}}
    \put(12,16){\makebox(0,0){\small $\gamma_0$}}
  \put(50,16){\makebox(0,0){\small $\gamma_2$}}
  \put(70,16){\makebox(0,0){\small $\gamma_3$}}
  \put(110,16){\makebox(0,0){\small $\gamma_{s-1}$}}
  \put(128,16){\makebox(0,0){\small $\gamma_s$}}

  \end{picture}
\end{center}
\bigskip

By construction the points of intersection of the curve $\gamma(w)$ with arcs in $\Gamma$ are indexed by the vertices of the string $w$.
The curve intersects the arcs of $\Gamma$ transversally, and since the string $w$ is reduced, none of the $\gamma_i$ is homotopic to a piece of an arc in $\Gamma$. 
Thus  the intersection numbers are minimal, and $I_{\Gamma}(\gamma(w)) = \ul{\dim}M(w)$.
Since $\gamma(w)$ has non-trivial intersection with arcs of $\Gamma$, it is clear that it is not homotopic to an arc in the triangulation $\Gamma$.
\bigskip

Conversely, let $\gamma: [0,1] \to S $ be a curve in $(S,M)$ which is not homotopic to an arc in $\Gamma$. We assume that the curve $\gamma$ is chosen (in its homotopy class) such that  it intersects the arcs $a$ of $\Gamma$ transversally (if at all) and such that the intersection numbers  $I_{\Gamma}(\gamma, a) $ are minimal.
Orienting $\gamma$ from $P=\gamma(0) \in M$ to $Q=\gamma(1) \in M$, we denote by $x_1$ the first internal arc of $\Gamma$ that intersects $\gamma$, by $x_2$ the second arc, and so on. We thus obtain a sequence $x_1, \ldots , x_s$ of (not necessarily different) internal arcs in $\Gamma$.
Since the intersection numbers are minimal, we know that $x_i \neq x_{i+1}$. Thus there are arrows, either 
$\alpha_i: x_i \to x_{i+1}$ or $\alpha_i: x_{i+1} \to x_{i}$ in $Q(\Gamma)$, and we obtain a walk
 $ w(\gamma) = \;x_1 \stackrel{\alpha_1}{\longleftrightarrow} x_2 \stackrel{\alpha_1}{\longleftrightarrow}\cdots  \stackrel{\alpha_{s-1}}{\longleftrightarrow} x_s$  in $Q(\Gamma)$.
The fact that $\gamma$ intersects the arcs of $\Gamma$ transversally implies that the walk $w(\gamma)$ is reduced and avoids the zero-relations, thus $w(\gamma)$ is a string in $A(\Gamma)$.

It follows from their construction that the two  maps  between strings and homotopy classes of curves defined above are mutually inverse.
\hfill $\Box$ \bigskip

{\bf Remark:}
Recall that two string modules $M(w)$ and $M(v)$ are isomorphic precisely when $v=w$ or $v=w^{-1}$. The inverse string $w^{-1}$ corresponds to orienting the curve in the opposite direction.

\begin{proposition}
Let $\Gamma$ be a triangulation of an unpunctured marked  surface $(S,M)$.  Then
there exists a bijection  between the homotopy classes of closed curves in $(S,M)$ and powers $b^n$ of bands  $b$ of $A(\Gamma)$.  
\end{proposition}

The proof is analogous to the proof  the previous  proposition.

\subsection{An example where $A(\Gamma)$ is not cluster-tilted}

We finally present in this section an example of an algebra $A(\Gamma)$ which is not cluster-tilted. Recall that an algebra $A$ is {\em tame} if for all $d \in \mathbb N$ there is a finite number $n_d$ of one-parameter families of $A-$modules such that almost every $d-$dimensional $A-$module belongs to one of these $n_d$ families. The algebra $A$ is said to be {\em domestic} if there is a constant $c$ such that  $n_d \le c$ for all $d \in \mathbb N$. On the other hand, if the numbers $n_d$ grow faster than any polynomial, then the tame algebra  $A$ is said to be of {\em non-polynomial growth}.
\medskip

It is well-known that every string algebra $A$ is tame, and that the one-parameter families are given by the bands in $A$ (see \cite{BR}). In particular, all the algebras $A(\Gamma)$ are tame,
because they are gentle and thus string algebras. 
Moreover the tame
cluster-tilted algebras of the form $A(\Gamma)$  studied here are all domestic, in fact they are of type $\mathbb{A}$ or $\tilde{\mathbb{A}}$, thus we may assume above that $c=0$ or $c=1$.
We construct in this section an example of an algebra $A(\Gamma)$ which is of non-polynomial growth, and thus cannot be cluster-tilted.
To obtain this example, we consider  a sphere $S$ with three holes and choose one marked point in each boundary component. We fix the following triangulation $\Gamma$ of $(S,M):$

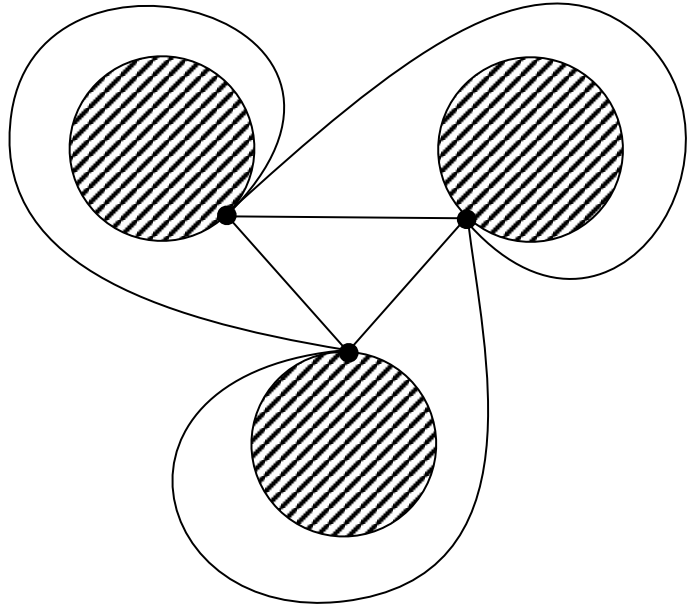
\begin{center}
\includegraphics[height=8cm]{sphere.pdf}
\end{center}

Then the algebra $A(\Gamma)$ is given by the following quiver with relations $\epsilon_i\rho_i=0, \rho_i\sigma_i=0$ and $\sigma_i\epsilon_i=0$ for $i = 1$ and $i=2$.
$$
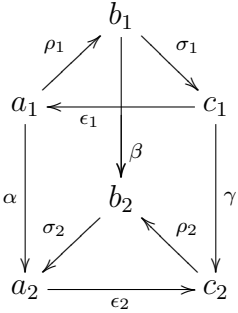
\xymatrix@R=20pt@C=20pt{
                                  & b_1 \ar[dd] \ar[dr]^{\sigma_1}  &                  \\
a_1 \ar[dd]_{\alpha} \ar[ur]^{\rho_1}    &                 \ar[d]^\beta               & c_1 \ar[ll]^{\epsilon_1\qquad}  \ar[dd]^{\gamma} \\
                                   & b_2 \ar[dl]_{\sigma_2}         &                                \\
a_2  \ar[rr]_{\epsilon_2}                &                                & c_2   \ar[ul]_{\rho_2}       \\
  }
$$

The string algebra $A(\Gamma)$ admits the following two bands 
$$ \xi = \;b_2 \stackrel{\sigma_2}{\longrightarrow} a_2 \stackrel{\alpha}{\longleftarrow} a_1  \stackrel{\rho_1}{\longrightarrow} b_1 \stackrel{\beta}{\longrightarrow} b_2$$
and
$$ \eta = \;b_2 \stackrel{\rho_2}{\longleftarrow} c_2 \stackrel{\gamma}{\longleftarrow} c_1  \stackrel{\sigma_1}{\longleftarrow} b_1 \stackrel{\beta}{\longrightarrow} b_2$$
\bigskip

Since $\xi$ and $\eta$ can be composed arbitrarily, the number of bands of a fixed length $l$ grows exponentially with $l$, and thus the algebra  $A(\Gamma)$ is of non-polynomial growth.
\bigskip

We would like to point out that the notion of non-polynomial growth of tame algebras discussed here does not coincide with the notion of non-polynomial growth cluster algebras discussed in \cite{FST}: There one counts the number of cluster variables, that is to say, of arcs in $(S,M)$, instead of one-parameter families, that is to say, closed curves in $(S,M)$. In  \cite{FST} the example we are considering in this section is classified as being of polynomial growth, meaning that, even if the number of curves is growing exponentially, the number of arcs is bounded for the sphere with three holes.

%
%

\bigskip

{\sc D\'epartement de
  Math\'ematiques,
  Universit\'e de Sherbrooke\\
  Sherbrooke (Qu\'ebec), J1K 2R1, Canada}  
  \bigskip

  {\sc Department of Mathematics, Bishop's University\\
    2600 College St., Sherbrooke, Quebec, Canada J1M 0C8}
\bigskip

{\sc  Universit\'e Denis Diderot (Paris VII)
Institut de MathŽmatiques\\
Case 7012-2, place Jussieu
75251 Paris Cedex 05}
\bigskip

  {\it  E-mail addresses: } \\
 {\tt ibrahim.assem@usherbrooke.ca}\\
 {\tt thomas.brustle@usherbrooke.ca} \quad and \quad
 {\tt tbruestl@ubishops.ca}\\
 {\tt gabrielle.charbonneau-jodoin@usherbrooke.ca}\\
{\tt plamondon@math.jussieu.fr}
 \medskip
\end{document}